%% file: PSIs.tex
\documentclass[review, 10pt]{article}         

\usepackage{mathtools}      
\usepackage[utf8]{inputenc} 
\usepackage[T1]{fontenc}    
\usepackage[hidelinks]{hyperref}       
\usepackage{url}            
\usepackage{microtype}      
\usepackage{amsfonts}       
\usepackage{nicefrac}       
\usepackage{amsthm}         
\usepackage{amssymb}        
\usepackage{bbm}            
\usepackage{booktabs}       
\usepackage{ctable}         
\usepackage{caption}        
\usepackage{multirow}       
\usepackage{tabularx,booktabs}  
\usepackage{graphicx}       
\usepackage{subcaption}
\usepackage{cleveref}       
\usepackage{cite}
\usepackage{array}          
\usepackage{dirtytalk}      
\usepackage{algorithm}          
\usepackage{algpseudocode} 
\usepackage{listings}
\usepackage{color}
\usepackage{array}
\usepackage[left=2.75cm,right=2.75cm, bottom=3.00cm, top=3.00cm]{geometry}
\usepackage{soul}
\usepackage[parfill]{parskip}
\usepackage{lineno}
\usepackage{setspace}
\usepackage{rotating}
\usepackage{subfiles}
\usepackage{enumitem}

\newcolumntype{Y}{>{\centering\arraybackslash}X}

\newcommand{\Tau}{\mathcal{T}} 

\newcommand{\norm}[1]{\left\lVert#1\right\rVert}

\newcommand\tab[1][0.5cm]{\hspace*{#1}}         

\providecommand{\keywords}[1]{
  \textbf{\textit{Keywords---}} #1
}

\newtheorem{definition}{Definition}[section]
\newtheorem{theorem}{Theorem}[section]
\newtheorem{corollary}{Corollary}[section]
\newtheorem{proposition}{Proposition}[section]

\theoremstyle{definition} 

\title{\bf Phase-space iterative solvers}
\author{Gaëtan Cortes,
Nur Cristian Sangiorgio,\\
Joaquin Garcia-Suarez\footnote{Corresponding author: \texttt{joaquin.garciasuarez@epfl.ch}}
}
\date{Institute of Civil Engineering,
\\ \'Ecole Polytechnique Fédérale de Lausanne (EPFL), CH 1015 Lausanne, Switzerland 
}           

\begin{document}

\maketitle

\begin{abstract}
We introduce an iterative method to solve problems in small-strain non-linear elasticity, termed ``Phase-Space Iterations'' (PSIs). 
The method is inspired by recent work in data-driven computational mechanics, 
which reformulated the classic boundary value problem of continuum mechanics using the concept of ``phase space'' associated with a mesh. 
The latter is an abstract metric space, whose coordinates are indexed by strains and stress components, where each possible state of the discretized body corresponds to a point. 
Since the phase space is associated to the discretized body, it is finite dimensional. 
Two subsets are then defined: an affine space termed ``physically-admissible set'' made up by those points that satisfy equilibrium and a ``materially-admissible set'' containing points that satisfy the constitutive law. 
Solving the boundary-value problem amounts to finding the intersection between these two subdomains. 
In the linear-elastic setting, this can be achieved through the solution of a set of linear equations; when material non-linearity enters the picture, such is not the case anymore and iterative solution approaches are necessary. 
Our iterative method consists on projecting points alternatively from one set to the other, until convergence. 
To evaluate the performance of the method, we draw inspiration from the ``method of alternative projections'' and the ``method of projections onto convex sets'', both of which have a robust mathematical foundation in terms of conditions for the existence of solutions and guarantees convergence. 
This foundation is leveraged to analyze the simplest case and to establish a geometric convergence rate. 
We also present a realistic case to illustrate PSIs' strengths when compared to the classic Newton-Raphson method, 
the usual tool of choice in non-linear continuum mechanics. 
Finally, its aptitude to deal with constitutive laws based on neural network is also showcased. 
\end{abstract}


\section*{List of symbols (in no particular order)}

\bgroup
\def\arraystretch{1.5}
\begin{table}[H]
\centering
\begin{tabularx}{0.95\textwidth}{ c *{3}{Y} }
& \textbf{Symbol} & \textbf{Meaning} &\textbf{Comment}  \\
\hline
& $n_e$ & Number of strain/stress components per element& \\
& $n_{\text{dofs}}$ & Number of nodal degrees of freedom components per node& \\
& $N_n$ & Number of nodes & \\
& $N_e$ & Number of elements & \\
& $N_{\text{dofs}}$ & Number of degrees of freedom & $ = N_{n} \cdot n_{\text{dofs}}$ 
\end{tabularx}
\end{table}
\egroup

\bgroup
\def\arraystretch{1.5}
\begin{table}[H]
\label{tab:glossary}
\centering
\begin{tabularx}{0.95\textwidth}{ c *{3}{Y} }
& \textbf{Symbol} & \textbf{Meaning} &\textbf{Comment}  \\
\hline
& $\boldsymbol{B}_e$ & Discrete gradient operator& $\in \mathbb{R}^{n_e \times N_{\text{dofs}}}$\\
& $\boldsymbol{B}_e^{\top}$ & Discrete divergence operator& $\in \mathbb{R}^{ N_{\text{dofs}} \times n_e }$\\
& $w_e$ & e-th element volume& \\
& $\boldsymbol{\sigma}_e$ & Physically-admissible stress (e-th element)& $\in \mathbb{R}^{n_e}$\\
& $\boldsymbol{\varepsilon}_e$ & Physically-admissible strain (e-th element)& $\in \mathbb{R}^{n_e}$\\
& $\boldsymbol{\sigma}'_e$ & Materially-admissible stress (e-th element)& $\in \mathbb{R}^{n_e}$\\
& $\boldsymbol{\varepsilon}'_e$ & Materially-admissible strain (e-th element)& $\in \mathbb{R}^{n_e}$\\
& $\boldsymbol{m}$ & Constitutive-law function & $\boldsymbol{\sigma}'_e=\boldsymbol{m}(\boldsymbol{\varepsilon}'_e)$\\
& $E$ & Physically-admissible set & Equilibrium and compatibility \\
& $P_E$ & Projection onto $E$ & \\
& $D$ & Materially-admissible set & Constitutive law \\
& $P_D$ & Projection onto $D$ & \\
& $Z_e$ & Local phase space & $\simeq \mathbb{R}^{2 n_e}$\\
& $Z$ & Global phase space & $\simeq \mathbb{R}^{2 n_e N_e}$\\
& $\boldsymbol{z}_e$ & Point in local phase space & $= \{ \boldsymbol{\sigma}_e,\boldsymbol{\varepsilon}_e \}\in Z_e$\\
& $\boldsymbol{z}$ & Point in global phase space & $= \{ [\boldsymbol{z}_e]_{e = 1}^{N_e} \}\in Z$\\
& $\boldsymbol{F}_{int}$ & Internal force vector & $\in \mathbb{R}^{N_{\text{dofs}}}$, it depends on stresses\\
& $\boldsymbol{F}_{ext}$ & External force vector & $\in \mathbb{R}^{N_{\text{dofs}}}$\\
& $\Delta \boldsymbol{F}$ & Force residual & $=\boldsymbol{F}_{int} - \boldsymbol{F}_{ext}$\\
& $\boldsymbol{C}$ & Distance constant matrix& \\
& $\boldsymbol{K}$ & Zero-strain stiffness matrix& \\
& $\mathbb{D}$ & Zero-strain moduli matrix& \\
& $Y$ & Young modulus& Function of strain\\
& $Y_0$ & Zero-strain Young modulus& \\
& $\nu$ & Poisson's ratio& \\
& $p$ & Non-linearity parameter& The larger, the more non-linear the material\\
& $\mathcal{X}, \mathcal{Y}, \mathcal{H}$ & finite-dimension Hilbert spaces&\\
& $T$ & a map between Hilbert spaces& \\
& $k$ & small parameter& dimensionless 
\end{tabularx}
\end{table}
\egroup

\section*{List of acronyms}

\bgroup
\def\arraystretch{1.5}
\begin{table}[H]
\centering
\begin{tabularx}{0.95\textwidth}{ c *{3}{Y} }
& \textbf{Acronym} & \textbf{Meaning} &\textbf{Comment}  \\
\hline
& DDCM & Data-driven computational mechanics& \\
& E-L & Euler-Lagrange & equations \\
& MAP & Method of alternating projections& \\
& NR & Newton-Raphson & method \\
& POCS & Projections onto convex sets & 
\\
& PSI & Phase-space iteration & 
\end{tabularx}
\end{table}
\egroup

\section{Introduction}
\label{Sec:Introduction}
\subfile{subfiles/Introduction}

\section{Method formulation}
\label{sec:formulation}
\subfile{subfiles/Method_formulation}

\section{Preliminary analysis of existence and uniqueness of solutions, and convergence properties of the method}
\label{sec:mathematical_props}
\subfile{subfiles/Uniqueness}


\subfile{subfiles/1D_bar}
\label{sec:1D_closed_form}

\section{Illustrative example: Kirchdoerfer's truss}
\label{sec:truss}
\subfile{subfiles/Kirchdoerfer}

\section{PSI combined with neural-network constitutive models}
\label{sec:NN_gaetan}
\subfile{subfiles/NeuralNetworks}

\section{Discussion}
\label{sec:discussion}
\subfile{subfiles/Discussion}

\section{Final remarks}
\label{sec:final}
\subfile{subfiles/Remarks_Acknow}

\bibliography{bibliography}{}
\bibliographystyle{unsrt}

\appendix

\section{Basic mathematical definitions}
\subfile{subfiles/appendix/Maths}


\section{Newton-Raphson (NR) method}
\subfile{subfiles/appendix/NR}
\label{app:NR}

\section{Neural-Network training}
\label{app:NeuralNetwork}
\subfile{subfiles/appendix/NeuralNetwork}

\end{document}

%% file: subfiles/Introduction.tex

Most practical problems in continuum mechanics lack a simple closed-form solution. 
Since the advent of modern computers, numerical methods have been developed to overcome this fact and to provide approximate solutions of particular problems for specific values of the relevant parameters \cite{bathe:2006}. 
The Newton-Raphson method (NR) has been a main player, due to its simplicity, ease of implementation and quadratic convergence \cite{Isaacson:1994,Oden:2013} 
(see \Cref{app:NR}). 
However, the method is not free from limitations. 
For once, it requires recomputing the stiffness of the system at every deformation level, what in practical terms means that ``tangent stiffness'' matrices have to be reassembled at each iteration. 
This can be time-consuming, and can also entail numerical issues if there are elements in the mesh whose lower stiffness worsens the conditioning of the overall matrix. 

A number of alternatives to NR have been proposed. 
Most of them are predicated on the minimization of the ``force residual'', i.e., the imbalance between internal and external forces, just like NR does, while simultaneously trying to avoid some of its limitations. 
For instance, quasi-Newton methods approximate successive inverse tangent matrices using the zero-strain stiffness and rank-one corrections by means of the Sherman-Morrison formula \cite{Zienkiewicz:1991}. 
Other methods iteratively scan the residual function locally, and search for the optimal step to minimize it. 
These are termed ``search methods'' \cite{Oden:2013,Zienkiewicz:1991}, the conjugate gradient method being one of the most popular flavors. 
Yet another approach, ``gradient flow'' (a.k.a. ``dynamic relaxation'') methods transform the elliptic problem of statics into an auxiliary parabolic problem \cite{Zienkiewicz:1991}. 
This can be solved with time-marching iterations \cite{Cook:2001} until a steady-state solution is achieved, and which is then taken to be the solution of the original problem. 

We introduce a new family of iterative solvers for boundary value problems in linear elasticity. 
The concept of ``phase space'' introduced in Ref.\,\cite{Trent_1} in the context of solvers for data-driven computational mechanics (DDCM) is borrowed. 
The phase space is $Z = \mathbb{R}^{2 n_e N_e}$, where $N_e$ is the number of elements in the mesh and $n_e$ is the number of independent components in either the strain or stress tensor in each element. 
Hence, each coordinate in the phase space corresponds to either a stress or a strain component in an element. 
To lighten this presentation, we are implicitly assuming one integration point per element (so stresses and strains are computed at a single location), and that all elements in the mesh are similar (e.g., only 1D bar elements). 
The global phase space $Z$ can be expressed as the Cartesian product of ``local'' phase spaces $Z_e=\mathbb{R}^{2 n_e}$ ($e=1, \ldots, N_e$), defined at the element level: $Z = Z_1 \times \ldots \times Z_{N_e}$.

Inspired by this seminal work of Kirchdoerfer and Ortiz in DDCM, 
we regard the exact solution of the problem as that that satisfies simultaneously boundary conditions, physical balance constraints (in our case Newton's second law, which boils down to static equilibrium in this case), 
kinematic considerations (the relation between strains and displacements), 
and a material-dependent relation between strain and stress. 
The former two combined defined an affine domain in the \textit{global} phase space that we term $E$, the ``physically-admissible set''. 
In those authors' original work, the ``materially-admissible set'' was known as a discrete set of points in the \textit{local} phase space. 
In this text, we will assume that, in each element, there is a function $\boldsymbol{m} : \mathbb{R}^{n_e} \to \mathbb{R}^{n_e}$ that defines a constitutive law $\boldsymbol{\sigma}_e' = \boldsymbol{m} (\boldsymbol{\varepsilon}_e')$. 
We term our method, inspired by DDCM, ``phase-space iterations'' (PSIs), which include DDCM as a particular case in which the constitutive information is known only at discrete points.

In Section 3, we introduce two mathematical theories with the potential to provide a general characterization of the method's performance. 
For this analysis, however, we have incorporate a pivotal simplification to be able to present an example that can be solved in closed form. This simplifcation enables us to apply these theories without requiring adaptations or introducing new concepts.
Specifically, we simplify by setting $N_e = n_e = 1$, so that $\{ \boldsymbol{\varepsilon}_e', \boldsymbol{m} (\boldsymbol{\varepsilon}_e') \} \subset Z_e$ defines a curve in $Z = Z_e = \mathbb{R}^2$, which, in turn, determines the limit points of a convex set.
The first theory we draw from is the \textit{method of alternating projection} (MAP), which dates back at least to Von Neumann's contributions \cite{Neumann:1949,Neumann:1951}  
(his iterative method was used for solving large systems of equations \cite{KACZMARZ:1993}, and has also been shown to have a strong connection to the method of subspace corrections \cite{Xu:1992,Xu:2002}). 
A similar philosophy was later applied to solve a related problem: finding the intersection of convex sets. 
The other theory is the \textit{method of projections onto convex sets} (POCS), which has an extensive literature \cite{Deutsch:1992}, with applications ranging from general inverse analysis \cite{Potter:1993} to image restoration \cite{Youla:1978,Demoment:1989}.
Convexity played a crucial role here, as it inherently ensured the desired mathematical properties, such as rapid convergence. 

NR is the alternative that the PSI method would have to contend with, so let us advance the advantages with respect to it:
\begin{enumerate}[leftmargin=*]
    \item Unlike Newton methods, PSI does not require assembling tangent stiffness matrices (or its inverse) at every step; an auxiliary matrix is assembled once, before any iteration, and it is reused later (see $\boldsymbol{K}$ in \Cref{sec:projection_onto_E}). 
    %
    %
    \item 
    %
    PSI is decomposed into two steps, one that involves solving a linear system of equations arising at the \textit{global} mesh level (its complexity scales with $N_e$), and another one that involves elementwise solution of a \textit{local} non-linear optimization (its complexity scales with $n_e$). 
    Half of the algorithm, i.e., the non-linear optimization step, is trivially parallelizable, meaning that each processor works on a completely independent task.  
    Conversely, 
    solvers like NR that act at the structure level requires ``domain decomposition'' \cite{FARHAT:1988,Hamandi:1995} to distribute the work between processors, what requires careful inter-processor communication \cite{MPI1994,open_mpi_1,open_mpi_2,open_mpi_3}.
    %
    %
    It is true that NR can still benefit from parallelization in other ways, e.g., using parallel linear solvers like MUMPS \cite{MPI1994,MUMPS:1,MUMPS:2} to solve the linear system of equations, but we note that this technique could also be used to solve the other (global) half of PSI. 
    \item Our distance-minimization method is more general, as it allows handling constitutive laws with less continuity. 
    NR requires derivatives of the constitutive law to be well-defined. 
    On the contrary, we can solve the distance minimization problem with derivative-less techniques (e.g., Brent's  method \cite{Brent:2013}). 
    %
\end{enumerate}

The paper is structured as follows: 
\Cref{sec:formulation} presents the method and simplifies the expressions for the particular case of a mesh made up by bar elements. 
\Cref{sec:mathematical_props} discusses existence and uniqueness of solutions, as well as the method's convergence rate, using simplified setting 
which can also be solved in closed form and analyzed in light of the theory. 
This section provides intuition as to the method performance and potential pitfalls. 
Generalization perspectives beyond the simplified system are briefly discussed, too. 
\Cref{sec:truss} applies the method to a complex truss case, and serves to parse out the influence of the parameters. 
Given the surging interest in constitutive laws based on neural networks \cite{review:2024}, \Cref{sec:NN_gaetan} explores how they fare in tandem with the new solver. 
A general discussion in \Cref{sec:discussion}, motivated by all the evidence gathered to that point, precedes \Cref{sec:final} which summarizes the paper and outlines future work directions.

\begin{figure}[H]
\centering
\captionsetup[subfigure]{justification=centering}
\begin{subfigure}[t]{.75\linewidth}
\caption{{\Large (a)} }
\includegraphics[width=\linewidth]{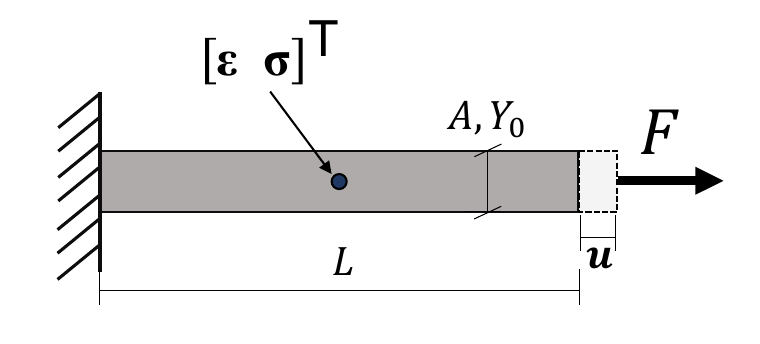}
\label{fig:scheme_1_bar}
\end{subfigure}

\begin{subfigure}[t]{.45\linewidth}
\caption{{\Large (b)} }
\includegraphics[width=\linewidth]{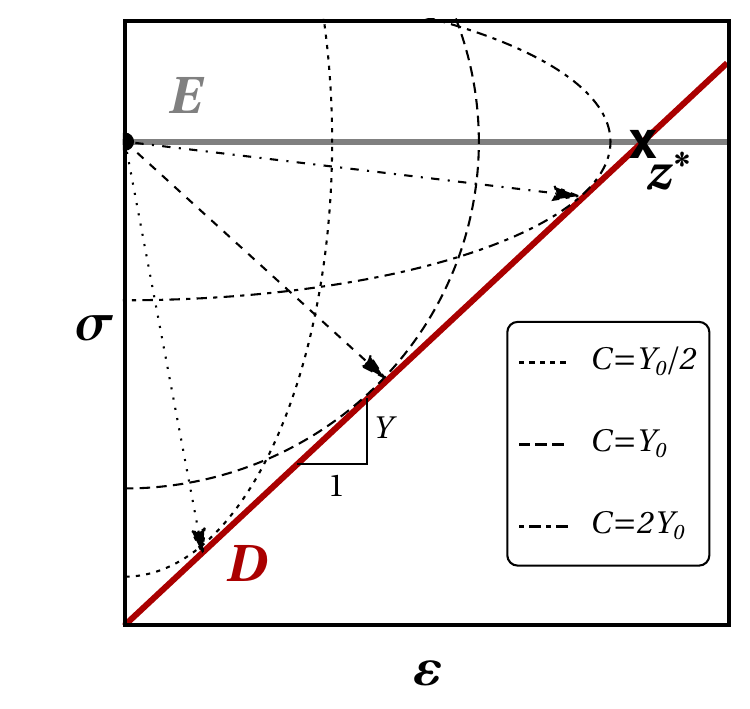}
\label{fig:iterations_PS}
\end{subfigure}
\hfill
\begin{subfigure}[t]{.45\linewidth}
\caption{{\Large (c)} }
\includegraphics[width=\linewidth]{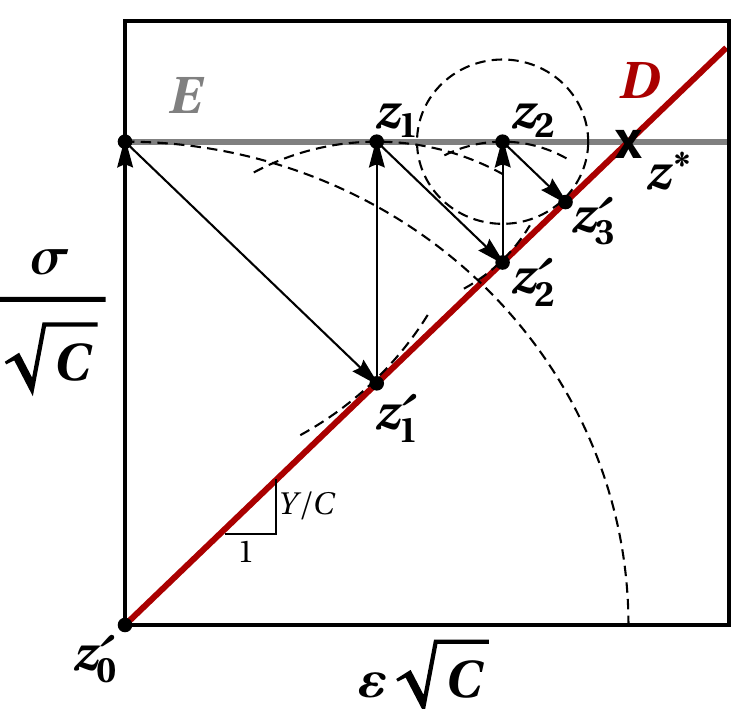}
\label{fig:iterations_PS_normalized}
\end{subfigure}
        \caption{Illustrating PSI method's key ideas. (a) Scheme of a one-element system whose phase-space is the simplest possible (will be used in \Cref{sec:1D_closed_form}). 
        (b) $E$ represents the physically-admissible set (equilibrium and compatibility), $D$ the materially-admissible one (constitutive law). 
        Phase-space possible projections: from a point satisfying equilibrium but not the constitutive law, $[0 \quad F/A ]^\top$, to different material-admissible points dependening on  choice of distance, i.e., parameter $C$ (see \cref{eq:distance_general,eq:distance} for details).  
        (c) Phase-space iterations in normalized space: every search/projection is equivalent to finding the point where the smallest possible circle centered in the previous iteration is tangent to the other set (only the first and the last circles are fully shown for clarity).}
        \label{fig:first}
\end{figure}

%% file: subfiles/Method_formulation.tex
Herein we present the new numerical method, leaving a brief presentation of NR to \Cref{app:NR}. 

\subsection{Introduction}

Consider a discretized body whose mesh is made up of $N_n$ nodes and $N_e$ elements, in which each node contains $n_{\text{dofs}}$ degrees of freedom (``dofs''), and stresses and strains are computed at one quadrature point (constant stress elements). 
In total, there are $N_{\text{dofs}} = n_{\text{dofs}} \times N_n$ degrees of freedom, some of them may be constrained while others are free or external forces are applied. 
Stresses and strains for each element are collected into vectors $\boldsymbol{\sigma}_e, \, \boldsymbol{\varepsilon}_e \in \mathbb{R}^{n_e}$, where $n_e$ is the minimal number of components to be considered ($n_e= 1$ for 1D elements, 3 for plane elasticity, 6 for 3D isotropic elasticity and 9 for 3D generalized continua \cite{Capriz:2013}). 
Thus, each element's local phase space $Z_e$ is isomorphic to $\mathbb{R}^{2 n_e}$ while the global phase space $Z$ is to $\mathbb{R}^{2 N_e n_e}$.

We then define a norm for elements in both the global ($\norm{\cdot}$) and local ($\norm{\cdot}_e$) phase space:
\begin{align}
    \norm{\boldsymbol{z}}
    =
    \left(
    \sum_{e = 1}^{N_e}
    w_e 
    \norm{\boldsymbol{z}}_e^2
    \right)^{1/2}
    =
    \left(
    \sum_{e = 1}^{N_e}
    {w_e \over 2} 
    \left[
    \boldsymbol{\varepsilon}_e^{\top}
    \boldsymbol{C}
    \boldsymbol{\varepsilon}_e
    +
    \boldsymbol{\sigma}_e^{\top}
    \boldsymbol{C}^{-1}
    \boldsymbol{\sigma}_e
    \right]
    \right)^{1/2}
    \, ,
    \label{eq:norm_general}
\end{align}
which naturally induces a metric and a distance (globally, $d$, and locally, $d_e$)
\begin{align}
    d (\boldsymbol{z}, \boldsymbol{z}')
    &=
    \left(
    \sum_{e = 1}^{N_e}
    w_e 
    d_e(\boldsymbol{z}_e,\boldsymbol{z}'_e)^2
    \right)^{1/2}
    =
    \norm{ \boldsymbol{z} - \boldsymbol{z}' }
    \nonumber \\
    &=
    \left(
    \sum_{e=1}^{N_e}
    {w_e \over 2} 
    \left[
    (\boldsymbol{\varepsilon}_e
    -
    \boldsymbol{\varepsilon}_e'
    )^{\top}
    \boldsymbol{C}
    (\boldsymbol{\varepsilon}_e
    -
    \boldsymbol{\varepsilon}_e'
    )
    +
    (\boldsymbol{\sigma}_e
    -
    \boldsymbol{\sigma}_e'
    )^{\top}
    \boldsymbol{C}^{-1}
    (\boldsymbol{\sigma}_e
    -
    \boldsymbol{\sigma}_e'
    )
    \right]
    \right)^{1/2} \, ,
    \label{eq:distance_general}
\end{align}
which resembles a Mahanalobis distance \cite{Maha:2002}, where positive definite $\boldsymbol{C}$ is a matrix with numbers with proper units so that the two addends are congruent, which is chosen to be both invertible and symmetric. 
$\boldsymbol{C}$ can be selected by the user, but we could use the zero-strain material moduli for both the global space distance when projecting onto $E$ and 
for projecting onto $D$. 
Keeping the matrix the same for both projections simplifies the analysis of the method. 
Note how the magnitude of $\boldsymbol{C}$ ``weights'' the two components of the distance: for relatively small values, the stress components would dominate the overall distance magnitude, while the strain ones would for large values of $\boldsymbol{C}$; this is illustrated by \Cref{fig:iterations_PS}.    

\subsubsection{Physical admissibility}

The static equilibrium, in discretized fashion, can be written as \cite{Trent_1}
\begin{align}
    \boldsymbol{F}_{int} (\boldsymbol{\sigma})
    -
    \boldsymbol{F}_{ext}
    =
    \sum_{e = 1}^{N_e}
    w_e \boldsymbol{B}_e^{\top} 
    \boldsymbol{\sigma}_e
    -
    \boldsymbol{F}_{ext}
    = 0 \, ,
    \label{eq:equilibrium}
\end{align}
where $\boldsymbol{F}_{ext} \in \mathbb{R}^{N_{\text{dofs}}}$ is the nodal force vector (containing both external forces and reactions), 
$w_e$ and $\boldsymbol{\sigma}_e$ are the volume of each element and its stresses (written as in vector form, i.e., Voigt notation), 
respectively, while $\boldsymbol{B}_e^{\top} \in \mathbb{R}^{N_{\text{dofs}} \times n_e}$ is the discrete divergence operator. 

Likewise, kinematic compatibility between displacements and infinitesimal strains $\boldsymbol{\varepsilon}_e$ (also a vector) in discrete form is
\begin{align}
    \boldsymbol{\varepsilon}_e 
    =
    \boldsymbol{B}_e
    \boldsymbol{u} \, ,
    \label{eq:compatibility}
\end{align}
where $\boldsymbol{u}$ is the nodal displacement field of that element and $\boldsymbol{B}_e \in \mathbb{R}^{n_e\times N_{\text{dofs}}} $ is the discrete gradient.

Thus, the first projection will take an initial phase space point satisfying the constitutive law 
(i.e., $\boldsymbol{z}' \in D = \{ [(\boldsymbol{\varepsilon}_e', \boldsymbol{\sigma}_e')]_{e=1}^{N_e} : \boldsymbol{\sigma}_e' - \boldsymbol{m} (\boldsymbol{\varepsilon}_e') = 0 \quad \forall e = 1, \ldots , N_e  \}$) 
to the closest, in the sense of \cref{eq:distance_general}, that belongs in the physically-admissible set, i.e.,  $\boldsymbol{z} \in  E = \{ [(\boldsymbol{\varepsilon}_e', \boldsymbol{\sigma}_e')]_{e=1}^{N_e} : \cref{eq:equilibrium} \}$. 
This can be expressed as
\begin{align}
    \boldsymbol{z}
    =
    P_E(\boldsymbol{z}')
    =
    \mathrm{argmin} \, \Pi_E(\cdot, \boldsymbol{z}') \, ,
    \label{eq:P_E}
\end{align}
where the functional is
\begin{align}
    \Pi_E(\boldsymbol{z}, \boldsymbol{z}')
    =
    d^2(\boldsymbol{z},\boldsymbol{z}')
    -
    \boldsymbol{\eta } \cdot
    \left( 
    \sum_{e = 1}^{N_e}
    w_e \boldsymbol{B}_e^{\top} 
    \boldsymbol{\sigma}_e
    -
    \boldsymbol{F}_{ext}
    \right) \, ,
    \label{eq:functional_Pi_E_general}
\end{align}
whose stationarity condition, i.e, Euler-Lagrange (E-L) equations \cite{gelfand:1963}, verify minimal distance to $\boldsymbol{z}'$ subject to equilibrium.
The Lagrange multipliers $\boldsymbol{\eta}$ can be thought to represent virtual displacements. 
%

See that compatibility, \cref{eq:compatibility}, can be replaced right away wherever $\boldsymbol{\varepsilon}_e'$ appears, thus the functional $\Pi_E$ depends on the field $\boldsymbol{\sigma}_e$, $\boldsymbol{\eta}$ (Lagrange multipliers) and $\boldsymbol{u}$ (displacement field). 


\subsubsection{Material admissibility}

The second projection returns the newly obtained $\boldsymbol{z} \in E$ back to $D$. 
Unlike the previous one, which is performed at the ``structure'' level in the global phase space (involving $\boldsymbol{z} = \{ [\boldsymbol{z}_e]_{e=1}^{N_e} \}$ and $\boldsymbol{z}' = \{ [\boldsymbol{z}_e']_{e=1}^{N_e} \}$), 
the new projection is the collection of $N_e$ projections performed at the element level in  the local phase spaces.  
In similar fashion (stretching the notation slightly --adding the tilde $' $-- to use $\boldsymbol{z}'$ to denote a new point on $D$ instead of $E$), we introduce a new functional, $\Pi_{D_e}$, defined element-wise, such that   
\begin{align}
    \Pi_{D_e}(\boldsymbol{z}_e', \boldsymbol{z}_e)
    =
    d_{e}^2(\boldsymbol{z}'_e
    ,
    \boldsymbol{z}_e)
    -
    \boldsymbol{\Lambda}_e
    \cdot
    \left( 
     \boldsymbol{\sigma}_e'
     - 
     \boldsymbol{m}(\boldsymbol{\varepsilon}_e')
    \right) \, ,
    \label{eq:functional_Pi_D_general}
\end{align}
$\boldsymbol{\Lambda}_e$ in this case representing virtual strains. Then,  
\begin{align}
    \boldsymbol{z}'_e
    =
    \mathrm{argmin} \, \Pi_{D_e}(\cdot, \boldsymbol{z}_e) \, ,
\end{align}
which must be solved $N_e$ times, one per element, but each of those sub-problems is of much less complexity than the global one. 
Hence, the second projection is defined as 
\begin{align}
    \boldsymbol{z}'
    =
    \{ [\boldsymbol{z}_e']_{e=1}^{N_e} \}
    =
    P_D(\boldsymbol{z})
    =
    \{ [
    \mathrm{argmin} \, \Pi_{D_e}(\cdot, \boldsymbol{z}_e)
    ]_{e=1}^{N_e} \} \, . \label{eq:general_projection_onto_D}
\end{align}

\subsubsection{Consecutive iterations}
\label{sec:consecutive_iterations}

We envision a geometric iterative solver in which each iteration consists of two projections, see \Cref{fig:iterations_PS,fig:iterations_PS_normalized}. 
After the two projections, the point obtained at the prior iteration (say, $\boldsymbol{z}'^{(n)}$) yields $\boldsymbol{z}'^{(n+1)} = P_D P_E (\boldsymbol{z}'^{(n)})$. 
%
%
Our preliminary analysis (\Cref{sec:mathematical_props}) suggests that the method, under some conditions, converges in norm, 
in the sense that the distance $\norm{(P_D P_E)^{n+1} \boldsymbol{z}'^{(0)} - (P_D P_E)^{n} \boldsymbol{z}'^{(0)}} \to 0 $ as $n \to \infty$. 

Unlike those mentioned in the introduction, see that the convergence of the method in this manner does not automatically imply that the force residual (i.e., the difference between external and internal forces) goes to zero in the same way. 
This is a remarkable difference when comparing to most solid mechanics solvers, which tend to be predicated in the minimization of the latter. 
For this reason, it is logical that PSI solvers may be equipped with a dual stop condition, which simultaneously checks both ``phase-space convergence'' (in terms of phase-space distance between iterations becoming smaller) and ``equilibrium convergence'' (in terms of force residual). Hence,  
the solution procedure is as follows: 
\begin{enumerate}[leftmargin=*]
    \item Choose a point $\boldsymbol{z}'^{(0)} \in Z$ that satisfies the constitutive laws defined at each element level, or simplify starting from the origin (which is part of the constitutive law).
    \item Apply $\boldsymbol{z}'^{(n)}
    =
    (P_D P_E)^{n} (\boldsymbol{z}'^{(0)})$ until either
    \begin{itemize}
        \item (a) the force residual is small enough, $\norm{\boldsymbol{F}_{int} (\boldsymbol{\sigma}') -\boldsymbol{F}_{ext}}_{L^2} < \mathrm{tol}_1 \norm{\boldsymbol{F}_{ext}}_{L^2}$, or
        \item (b) the relative phase-space distance between consecutive iterations is small enough, \newline $\norm{\boldsymbol{z}'^{(n+1)} - \boldsymbol{z}'^{(n)}} < \mathrm{tol}_2 \norm{\boldsymbol{z}'^{(n)}}$,
    \end{itemize}
\end{enumerate}
$\mathrm{tol}_2 =  \mathrm{tol}_1 / 10$ proved satisfactory based on observations.

\subsection{Similarities with other methods}


Ref. \cite{Ladeveze:2012} considers a large class of iterative methods with ``two search directions'', particularly in the context of quasi-static time evolution, i.e., fields can evolve irreversibly in time, but inertia is neglected. 
A similar problem piecing to the one we are using is common to all members of this class: 
a global linear problem (equilibrium) and a non-linear local problem (constitutive law, where the strain rates can also appear), 
the final solution is found iterating from one to the other using the two search directions. 
The NR and PSI methods fall into this category, and so does the LATIN method (``LArge Time INcrements''), 
they differ in the way the two subdomains are defined and in the way the search directions are chosen. 
%
We remark that similarities with MAP and POCS were not leveraged.  

\subsection{Operative form of the projections}

\subsubsection{Projection onto $E$}
\label{sec:projection_onto_E}

%
%
Assume $\boldsymbol{z'} = \{ (\boldsymbol{\varepsilon}_e',\boldsymbol{\sigma}_e')_{e=1}^{N_e} \}$ is given. 
Enforcing the stationarity condition ($\delta \Pi_E = 0$) to obtain, field by field,
\begin{subequations}
\begin{align}
    \delta \boldsymbol{u} & 
    \to 
    \sum_{e=1}^{N_e}
    w_e
    \boldsymbol{B}_{e}^{\top}
    \boldsymbol{C}
    (
    \boldsymbol{B}_{e} \boldsymbol{u}
    -
    \boldsymbol{\varepsilon}_e'    
    )
    = 0 
    \to
    \left(
    \sum_{e = 1}^{N_e}
    w_e 
    \boldsymbol{B}_{e}^{\top}
    \boldsymbol{C}
    \boldsymbol{B}_{e} 
    \right)
    \boldsymbol{u}
    =
    \sum_{e = 1}^{N_e}
    w_e 
    \boldsymbol{B}_{e}^{\top}
    \boldsymbol{C}
    \boldsymbol{\varepsilon}_e'
    \, , \label{eq:stationarity_u_k_general}
    \\
    \delta \boldsymbol{\sigma}_e & \to 
    w_e
    \boldsymbol{C}^{-1}
    \left(
    \boldsymbol{\sigma}_e - \boldsymbol{\sigma}_e'
    \right)
    -  
    w_e 
    \boldsymbol{B}_{e} \boldsymbol{\eta}
    = 0 
    \to 
    \boldsymbol{\sigma}_e 
    =
    \boldsymbol{\sigma}_e'
    +
    \boldsymbol{C}
    \boldsymbol{B}_{e} \boldsymbol{\eta}
    \, , \label{eq:stationarity_sigma_general}
    \\
    \delta \boldsymbol{\eta} & 
    \to
    \sum_{e = 1}^{N_e}
    w_e 
    \boldsymbol{B}_{e}^{\top} 
    \sigma_{e}
    -
    \boldsymbol{F}_{ext}
    = 0
    \, . \label{eq:stationarity_eta_general}
\end{align}
\label{eq:stationarity}
\end{subequations}    

Upon combination of \cref{eq:stationarity_sigma_general} with \cref{eq:stationarity_eta_general}, one obtains:
\begin{align}
    \left(
    \sum_{e = 1}^{N_e}
    w_e 
    \boldsymbol{B}_{e}^{\top}
    \boldsymbol{C}
    \boldsymbol{B}_{e} 
    \right)
    \boldsymbol{\eta}
    =
    \boldsymbol{K}
    \boldsymbol{\eta}
    =
    \boldsymbol{F}_{ext}
    -
    \sum_{e = 1}^{N_e}
    w_e 
    \boldsymbol{B}_{e}^{\top}
    \sigma_e' 
    =
    \boldsymbol{F}_{ext}
    -
    \boldsymbol{F}_{int} \nonumber \\
    \to
    \boldsymbol{\eta}
    =
    \boldsymbol{K}^{-1}
    \left(
    \boldsymbol{F}_{ext}
    -
    \boldsymbol{F}_{int}(\boldsymbol{\sigma}')
    \right)
    =
    \boldsymbol{K}^{-1}
    \Delta \boldsymbol{F}(\boldsymbol{\sigma}')
    \, , \label{eq:sol_eta}
\end{align}
wherefrom it becomes apparent that $\boldsymbol{\eta}$ can be understood as a measure of the imbalance between external forces and internal forces \textit{for a given value of materially-admissible stresses} $\boldsymbol{\sigma}'$. 
In practice, the equation is solved only for those degrees of freedom that come not imposed by boundary conditions. For the latter, it is set from the start that $\eta_i = 0$.
Once the nodal imbalance variables $\boldsymbol{\eta}$ are available, from \cref{eq:stationarity_sigma_general} we find
\begin{align}
    \boldsymbol{\sigma}_e
    =
    \boldsymbol{\sigma}_e'
    +
    \boldsymbol{C}
    \boldsymbol{B}_{e} 
    \boldsymbol{K}^{-1}
    \left(
    \boldsymbol{F}_{ext}
    -
    \boldsymbol{F}_{int}(\boldsymbol{\sigma}')
    \right) \, , 
    \label{eq:sigma_update_P_E}
\end{align}
for each element. 

The essential boundary conditions must be enforced as part of finding the physically-admissible solution. 
We can set another system of linear equations from \cref{eq:stationarity_u_k_general}, thus
\begin{align}
    \boldsymbol{u}
     =
     \boldsymbol{K}^{-1}
     \left(
     \sum_{e=1}^{N_e}
    w_e
    \boldsymbol{B}_{e}^{\top}
    \boldsymbol{C}
    \boldsymbol{\varepsilon}'    
    \right)
     \, ,   
     \label{eq:sol_u}
\end{align}
and then enforce the essential BCs ($u_i = \hat{u}_i$) for some prescribed value $\hat{u}_i$ . 
In the implementation of the code used for \Cref{sec:truss,sec:NN_gaetan}, this is done by ``condensing'' the imposed displacements, 
substituting them directly into the vector $\boldsymbol{u}$ and solving \cref{eq:sol_u} for a reduced system that includes the ``free'' degrees of freedom and forces arising from condensation \cite{bathe:2006,Cook:2001}.

%
%
%

Summarizing, the first projection goes from $\boldsymbol{z}' 
= 
\{ [(\boldsymbol{\varepsilon}'_e,\boldsymbol{\sigma}_e')]_{e = 1}^{N_e} \} \in D$ to 
\begin{align}
    \boldsymbol{z}
    &= 
    \{ [(\boldsymbol{\varepsilon}_e,\boldsymbol{\sigma}_e)]_{e = 1}^{N_e} \} = 
    P_E (\boldsymbol{z}') 
    = 
    \{ [
    \boldsymbol{B}_{e}
    \boldsymbol{K}^{-1}
     \sum_{e=1}^{N_e}
    w_e
    \boldsymbol{B}_{e}^{\top}
    \boldsymbol{C}
    \boldsymbol{\varepsilon}'
    ,
    \boldsymbol{\sigma}_e' 
    +
    \boldsymbol{C}
    \boldsymbol{B}_{e} 
    \boldsymbol{K}^{-1}
    \Delta \boldsymbol{F}(\boldsymbol{\sigma}')]_{e = 1}^{N_e} \}
    \in E \, .
\end{align}
%


See that the matrix $\boldsymbol{K}$ depends on the choice of $\boldsymbol{C}$. 
For instance, when using the zero-strain elastic constants, $\boldsymbol{K}$ becomes the tangent matrix at the origin, which is computed and stored once and needs no updating across iterations. 


\subsubsection{Projection onto $D$}

Assume $\boldsymbol{z} = \{ (\boldsymbol{\varepsilon}_e,\boldsymbol{\sigma}_e)_{e=1}^{N_e} \} ] \in E$ is given. 
When it comes to the projection over the constitutive law, the stationarity condition for \cref{eq:functional_Pi_D_general} yields the following E-L equations 
\begin{subequations}
\begin{align}
    \delta \boldsymbol{\varepsilon}_e' & 
    \to 
    \boldsymbol{C}
    \left(
    \boldsymbol{\varepsilon}_e' - \boldsymbol{\varepsilon}_e
    \right)
    +  
    \boldsymbol{\Lambda}_e
    \nabla \boldsymbol{m} (\boldsymbol{\varepsilon}_e')
    = 0 
    \to
    \boldsymbol{\varepsilon}_e' 
    =
    \boldsymbol{\varepsilon}_e
    -
    \boldsymbol{C}^{-1}
    \boldsymbol{\Lambda}_e
    \nabla \boldsymbol{m} (\boldsymbol{\varepsilon}_e')
    \, , \label{eq:stationarity_eps_general_D}
    \\
    \delta \boldsymbol{\sigma}_e' & 
    \to 
    \boldsymbol{C}^{-1}
    \left(
    \boldsymbol{\sigma}_e' - \boldsymbol{\sigma}_e
    \right)
    -  
    \boldsymbol{\Lambda}_e
    = 0 
    \to 
    \boldsymbol{\Lambda}_e
    =
    \boldsymbol{C}^{-1}
    \left(
    \boldsymbol{\sigma}_e' - \boldsymbol{\sigma}_e
    \right)
    \, , \label{eq:stationarity_sigma_general_D}
    \\
    \delta \boldsymbol{\Lambda}_e & 
    \to
    \boldsymbol{\sigma}_e'
    -
    \boldsymbol{m}(\boldsymbol{\varepsilon}_e')
    = 0
    \to
    \boldsymbol{\sigma}_e'
    =
    \boldsymbol{m}(\boldsymbol{\varepsilon}_e')
    \, . \label{eq:stationarity_eta_general_D}
\end{align}
\end{subequations}    
Enough continuity has to be assumed so that the gradient of the constitutive law ($\nabla \boldsymbol{m}$) is well-defined. 
See that we can minimize the distance directly, through some numerical approach, 
after enforcing the constitutive law relation \cref{eq:stationarity_eta_general_D}. 
Minimizing this distance is precisely the traditional approach in DDCM \cite{Trent_1}, 
but in this case $\boldsymbol{m}(\varepsilon_e')$ is not known, so one has to scan a discrete dataset. In our case $\boldsymbol{m}(\varepsilon_e')$ is defined, so
functional minimization of the objective function, based on its derivatives or not, can be leveraged to find the minimizer directly. 

Combining \cref{eq:stationarity_eps_general_D,eq:stationarity_sigma_general_D,eq:stationarity_eta_general_D} yields a non-linear vector equation for $\boldsymbol{\varepsilon}_e'$: 
\begin{align}
    \boldsymbol{\varepsilon}_e' 
    =
    \boldsymbol{\varepsilon}_e
    +
    \boldsymbol{C}^{-1}
    \boldsymbol{C}^{-1}
    (\boldsymbol{\sigma}_e
    - \boldsymbol{m}(\boldsymbol{\varepsilon}_e'))
    \nabla \boldsymbol{m} (\boldsymbol{\varepsilon}_e')
    \, .
    \label{eq:E-L_optimal_C_general}
\end{align}
This equation is referred to as ``the Euler-Lagrange equation of $P_D$'' hereafter, as it combines \cref{eq:stationarity_eps_general_D,eq:stationarity_sigma_general_D,eq:stationarity_eta_general_D} into one. 
It represents a vector equation whose complexity depends on the form of the function $\boldsymbol{m}$, and 
on the number of components of the vectors ($n_e$), which depends on the problem: for instance, 3 components for plane stress or strain, and for 1D-element meshes it boils down to a scalar equation.

\subsection{Simplifications for 1D elements (bars)}




By virtue of the typology of the connection among bar elements, they work primarily by axially stretching or compressing, so the mechanical state of the $e$-th element corresponds to the simplest local phase space, i.e., $Z_e = \mathbb{R}^2$ ($n_e= 1$), each point defined by only two coordinates $(\sigma_e, \varepsilon_e)$. 
See that $\boldsymbol{C}$ becomes a scalar $C$. 
In this case, the norms and distances, \cref{eq:norm_general,eq:distance_general} simplify considerably:
\begin{align}
    \norm{\boldsymbol{z}}
    &=
    \left(
    \sum_{e = 1}^{N_e}
    w_e 
    \norm{\boldsymbol{z}}_e^2
    \right)^{1/2}
    =
    \left(
    \sum_{e = 1}^{N_e}
    w_e 
    \left[
    {\varepsilon_e^2 \over 2 C^{-1}}
    +
    {\sigma_e^2 \over 2 C}
    \right]
    \right)^{1/2}
    \, , \\
    d (\boldsymbol{z}, \boldsymbol{z}')
    &=
    \norm{ \boldsymbol{z} - \boldsymbol{z}' }
    =
    \left(
    \sum_{e=1}^{N_e}
    w_e
    \left[
    { (\sigma_e - \sigma_e')^2 
    \over 
    2 C}
    +
    { (\varepsilon_e - \varepsilon_e')^2 
    \over 
    2 C^{-1}}
    \right]
    \right)^{1/2} \, ,
    \label{eq:distance}
\end{align}
where, recall, $C$ is but a number with proper units so that the two addends in brackets are congruent. 
%
The projection onto $E$ does not change qualitatively,  
when it comes to $D$ (enforcing the constitutive law relation), we obtain:
\begin{align}
    \boldsymbol{z}_e'
    =
    \mathrm{arg min} 
    \left\{
    {(\varepsilon_e' - \varepsilon_e)^2 
    \over 
    2 C^{-1}}
    +
    {(m(\varepsilon_e') - \sigma_e)^2 
    \over 
    2 C}
    \right\} \, .
    \label{eq:reduced_functional_D}
\end{align}
and the Euler-Lagrange equation, similarly to \cref{eq:E-L_optimal_C_general}, yields
\begin{align}
    \varepsilon_e' = \varepsilon_e
    +
    {(\sigma_e - m(\varepsilon_e')) 
    \over 
    C^2}
    {d m(\varepsilon_e') \over d \varepsilon_e'}
    \, .
    \label{eq:E-L_C_truss}
\end{align}




\subsection{Algorithm pseudo-code}

\begin{algorithm}[H]
\caption{Phase-space iterative projections}
\begin{algorithmic}
\State \textbf{Require:} Connectivity, $\forall e = 1 , \ldots, N_e$, compatibility matrices $\boldsymbol{B}_e$; zero-strain moduli's matrix $\mathbb{D}$; matrix of distance constants $\boldsymbol{C}$; $\forall i = 1 , \ldots, N_{\text{dof}}$, external forces $F_{ext}$ or boundary conditions. 

\State (i) Set $k=0$. Initial data assignation:
\For{$  e = 1 , \ldots, N_e$ }\\
\tab Set $\boldsymbol{z'}_e^{(0)} =  (0, 0)$, unless there being essential boundary conditions in the nodes of the said element, in that case set $\boldsymbol{\varepsilon}_e' = \boldsymbol{B}_e \hat{\boldsymbol{u}}$. 
\EndFor

\State (ii) Project onto $\mathrm{E}$: \\
Given $\boldsymbol{z'}^{(n)}$, find $\boldsymbol{z}^{(n+1)}=\{ \boldsymbol{z}_e^{(n+1)} \}_{e=1}^{N_e} =\{ (\varepsilon_e,\sigma_e) \}_{e=1}^{N_e}$ from solving \cref{eq:stationarity} to get $\boldsymbol{u}^{(n+1)}$ and $\boldsymbol{\eta}^{(n+1)}$ 
\For{$ e = 1, \ldots , N_e$ }\\
\tab Compute $\boldsymbol{\varepsilon}_e$ from \cref{eq:compatibility} and  $\boldsymbol{\sigma}_e$ from \cref{eq:sigma_update_P_E}
\EndFor

\State (iii) Project onto $\mathrm{D}$: 
\For{$ e =1, \ldots, N_e$ }\\
\tab Given $\boldsymbol{z}^{(n+1)} = \{ \boldsymbol{z}^{(n+1)}_e \}_{e=1}^{N_e}$, find $\boldsymbol{z'}^{(n+1)}_e = (\boldsymbol{\varepsilon}'_e,\boldsymbol{\sigma}'_e)$, where $\boldsymbol{\varepsilon}'_e$ is got from \cref{eq:general_projection_onto_D} and then  $\boldsymbol{\sigma}'_e = \boldsymbol{m}(\boldsymbol{\varepsilon}'_e)$. 
\EndFor

\State (iv) Convergence test:
\If{ $\boldsymbol{K} \boldsymbol{\eta}^{(n+1)} < \mathrm{tol}_1 \norm{\boldsymbol{F}_{ext}}_{L^2}$} \\
\tab Final displacement field $\boldsymbol{u} = \boldsymbol{u}^{(n+1)}$ \textbf{exit}\\
\textbf{else if $\norm{\boldsymbol{z}^{(n+1)} - \boldsymbol{z}^{(n)}} < \mathrm{tol}_2 \norm{ \boldsymbol{z}^{(n)}} $} \\
\tab Final displacement field $\boldsymbol{u} = \boldsymbol{u}^{(n+1)}$ \textbf{exit}
\Else \\
\tab $n \leftarrow n+1$, \textbf{goto} (ii)
\EndIf
\end{algorithmic}
\label{alg:PSIs}
\end{algorithm}

%% file: subfiles/Uniqueness.tex
In this section, we provide a preliminary analysis of the method. 
Notions of existence and uniqueness are borrowed from POCS, while the convergence rate analysis is inspired by the MAP. 
It is just a preliminary one because, to leverage existing results from the literature without introducing additional modifications, the scope of this section is restricted to a single 1D bar element (length $L$, cross-section $A$ and Young's modulus $Y$), \Cref{fig:scheme_1_bar}. 
This approach is convenient because the example can be solved in closed form, allowing the abstract aspects to be directly compared with explicit computations. 
Hence, since there is a single element, local and global phase-space overlap and satisfy $Z_e = Z \simeq \mathbb{R}^2$. 
Stretching the notation, we will still use $D$ to refer to the convex set $D = \{ (\varepsilon', \sigma') \in Z : \sigma' - m(\varepsilon') \le 0 \}$, whose upper frontier is the constitutive law, and $E$ to refer to the convex set  $E = \{ (\varepsilon', \sigma') \in Z : \sigma' - F/A \ge 0 \}$, whose lower frontier is the set of points that satisfy equilibrium and compatibility. 
The discussion on how to generalize this simplest case to other scenarios is presented at the end of the section.

\subsection{Preliminaries}

Let us begin by presenting some relevant results for maps between Hilbert spaces, to then apply them to projections. 
Since the method is iterative, we review some important properties of maps that will be necessary, namely, \textit{non-expansiveness}, \textit{asymptotic regularity} and the set of \textit{fixed points} of a map. 
These concepts, along with the one of convex subset, are combined with Opial's theorem (its finite-dimensional version for our purposes). 
Proofs of theorems, propositions and corollaries that are not explicitly references in this section can be found in \cite{gubin1966method,lewismalick,MAP,convrate}. 

\vspace{0.25cm}

\begin{definition}[Non-expansiveness]
Given $\mathcal{X},\mathcal{Y}$ normed vector spaces, a map $T : \mathcal{X} \rightarrow \mathcal{Y}$ is \textit{non-expansive} if
\begin{align}
    \|T(\boldsymbol{x}_1)-T(\boldsymbol{x}_2)\| 
    \leq 
    \|\boldsymbol{x}_1-\boldsymbol{x}_2\|, \
    \forall 
    \boldsymbol{x}_1,\boldsymbol{x}_2 
    \in \mathcal{X} \, . 
\nonumber
\end{align}  
\end{definition} 

\vspace{0.25cm}

\begin{definition}[Asymptotic regularity]
Let $T^{(n)}(\boldsymbol{x})$ represent n successive compositions of $T$ applied to $\boldsymbol{x}$. 
A map $T$ of a metric space $(M,d)$ onto itself is said to be asymptotically regular if 
\begin{align}
    \lim_{n \to \infty} 
    d\left(
        T^{(n)}(\boldsymbol{x}), 
        T^{(n+1)}(\boldsymbol{x})
    \right) 
    = 0 \, , 
    \forall \boldsymbol{x}\in M
    \, .
    \nonumber
\end{align}
\end{definition}

\vspace{0.25cm}

\begin{definition}[Set of fixed points $\Tau_T$]
$\Tau_T = \{ \boldsymbol{x} \in \mathcal{X} : T(\boldsymbol{x})=\boldsymbol{x} \}$ is the set of fixed points of $T$. 
\end{definition}

\vspace{0.25cm}

\begin{theorem}[Opial]
Let $\mathcal{C}$ be a closed convex subset of a finite dimensional inner product space $\mathcal{X}$. 
Let $T : \mathcal{C} \rightarrow \mathcal{C}$ be any nonexpansive, asymptotycally-regular map. 
Let $\Tau_T$ denote the set of fixed points of $T$. 
Assume $\Tau \neq \emptyset$. 
Then, for every $\boldsymbol{x} \in \mathcal{C}$, the sequence $T^{(n)}(\boldsymbol{x})$ (n-th successive composition of $T$) converges to a point in $\Tau_T$.  
\label{Opial}
\end{theorem}

Next, we review the definition of projections onto convex sets so that we can relate the previous facts to projections later. 
\vspace{0.25cm}


\begin{theorem}[Hilbert projection theorem]
For every vector 
$\boldsymbol{x}$ in a Hilbert space $\mathcal{H}$ and every nonempty closed convex set $\mathcal{C} \subseteq \mathcal{H}$ there exists a unique vector $\hat{\boldsymbol{x}} \in \mathcal{C}$ for which $\|\boldsymbol{x} - \hat{\boldsymbol{x}}\|$ = $\displaystyle \inf _{c\in \mathcal{C}}\|\boldsymbol{x}-c\|$. 
\label{thm:HilbertProjection}
\end{theorem}

\vspace{0.25cm}

\begin{definition}[Projection Operator]
The projection operator onto a closed-convex subset $\mathcal{C}_i$ is defined as $P_i(\boldsymbol{x}):=\hat{\boldsymbol{x}}$ (see \Cref{thm:HilbertProjection})
\label{def:ProjectionOperator}
\end{definition}

\vspace{0.25cm}

\begin{definition}[Intersection of convex subsets $\bar{\mathcal{C}}$]
   Let $\mathcal{C}_1,\ldots,\mathcal{C}_m$ be closed convex subsets of $\mathcal{H}$, then $\bar{\mathcal{C}}:=\bigcap_{i=1}^{m} \mathcal{C}_i$ 
\end{definition}

\vspace{0.25cm}

\begin{definition}[Global projection $P(\boldsymbol{x})$]
Let $P_1,...,P_m$ be the projection operators corresponding to the closed convex subsets $\mathcal{C}_1,...,\mathcal{C}_m$, the overall projection is 
$P(\boldsymbol{x}):=P_m(P_{m-1}(...(P_2(P_1(\boldsymbol{x})))))$. 

\label{combinationprojections}
\end{definition}

\vspace{0.25cm}

\begin{proposition}
If $\mathcal{C}$ is a closed convex subset of a real Hilbert space $\mathcal{H}$, then, for every $\boldsymbol{x}_1, \boldsymbol{x}_2 \in \mathcal{H}: $ $\|\hat{\boldsymbol{x}_1}-\hat{\boldsymbol{x}_2}\|^2 \leq \langle \boldsymbol{x}_1-\boldsymbol{x}_2, \hat{\boldsymbol{x}_1}-\hat{\boldsymbol{x}_2} \rangle$, which implies also, by Cauchy-Schwarz inequality,
that projection operators are \textit{non-expansive}, $\|P(\boldsymbol{x}_1)-P(\boldsymbol{x}_2)\|\leq\|\boldsymbol{x}_1-\boldsymbol{x}_2\|$ \label{prop3}
\end{proposition}

\vspace{0.25cm}

\begin{proposition}
$P(\boldsymbol{x}):=P_m(P_{m-1}(...(P_2(P_1(\boldsymbol{x})))))$ is \textit{non-expansive}, since $P_i$ is \textit{non-expansive} $\forall i=1,\dots,m$ (from proposition \ref{prop3})
\label{combined-nonexpansiveness}
\end{proposition}

The final proposition and corollary wrap this review up, by showing how successive application of the projection operator, which itself can be a composition of $m$ atomic projections, converges (in the sense of \textit{asymptotic regularity}) to the intersection of the closed convex sets that define the atomic projections. 

\vspace{0.25cm}

\begin{proposition}
Assume $\bar{\mathcal{C}} \neq \emptyset$, $P(\boldsymbol{x})=P_m(P_{m-1}(...(P_2(P_1(\boldsymbol{x})))))$, then $\bar{\mathcal{C}}=\Tau_P$ (i.e., the fixed points of the projection are the intersection of the closed convex subsets). 
\label{fixedpointsandsolution}
\end{proposition}

\vspace{0.25cm}

\begin{corollary}
    Assume $\bar{\mathcal{C}} \neq \emptyset$, $P(\boldsymbol{x})=P_m(P_{m-1}(...(P_2(P_1(\boldsymbol{x})))))$. 
    It follows that $P(\boldsymbol{x})$ is \textit{asymptotically regular}.  
    \label{corollaryasymptoticregularity}
\end{corollary}

\subsection{Proof of existence of a fixed point and of convergence}

Here, we will show proof of existence of a solution and of convergence of the iterative approach to a fixed point.

\vspace{0.25cm}

\begin{theorem}[Convergence of the algorithm to a fixed point] 
Let the materially-admissible set $D$ and the physically-admissible set $E$ represent two closed convex sets in the phase space defined by (\ref{eq:norm_general}). 
Assume that their intersection $\bar{\mathcal{C}} = E \cap D \neq \emptyset$, 
and that in \cref{alg:PSIs}, the matrix $\boldsymbol{C}$ is symmetric and positive definite. 
Then, in the limit of infinite iterations, \cref{alg:PSIs} is guaranteed to converge to a fixed point in $\bar{\mathcal{C}}$. 
    
\end{theorem}
\begin{proof}
We will divide the proof in steps. \Cref{fig:iterations_PS,fig:iterations_PS_normalized} may serve as graphical illustration. 

\paragraph{Step 1: show that the phase space $Z$ is a Hilbert space}
This can be formally shown using an isomorphism argument (for any $N_e$, $n_e$) between the phase-space and $\mathbb{R}^{2 n_e N_e}$, which is obviously a finite dimensional Hilbert space. 
Indeed, $Z_e \simeq \mathbb{R}^{2 n_e}$ and $Z \simeq \mathbb{R}^{2 n_e N_e}$, 
and under some necessary assumptions on the matrix $\boldsymbol{C}$, the phase-space with respect to $\mathbb{R}^{2 n_e N_e}$ is a combination of a rescaling of axis and of a change of basis through an orthogonal matrix, using the well known spectral theorem for symmetric matrices \cite{Axler:2024}.

\paragraph{Step 2: show that $P_E$ and $P_D$ are projection operators}
In other words, show that they are consistent with \Cref{def:ProjectionOperator}:

\begin{itemize}[leftmargin=*]
    \item The minimization of the functionals given by eqs.(\ref{eq:functional_Pi_E_general}) and (\ref{eq:functional_Pi_D_general}) is equivalent to minimizing the distance between a given point in the phase space and any other point, under the constraint of belonging to sets $E, D$ respectively (operatively, those constraints are added via Lagrange multipliers), cf. \Cref{thm:HilbertProjection}. 
    \item Under the assumptions on the convexity of $E, D$, we can use theorem \ref{thm:HilbertProjection} and definition \ref{def:ProjectionOperator} in order to realize that the minimization of these functionals in \Cref{alg:PSIs} is exactly the definition of a projection operator onto a convex set if we understand those (semi-infinite) convex sets to be defined by limit points pertaining to the equilibrium set
    and 
    to the constutive law
    along with the proper axes,
    %
    see \Cref{fig:convex_sets}. 
    This fact may appear to limit the method as it stands right now, given that it is equivalent to assuming that the positive/negative stress components associated to positive/negative strains remain of the same sign across iterations, but this is a valid supposition, see \Cref{fig:evolution_projections} in \Cref{sec:influence_distance}.  
    %
    \item Defining the combination of the projection onto $D$ and $E$ as $P(\boldsymbol{x})$ in accordance with definition \ref{combinationprojections}, we also have that $P(\boldsymbol{x})$ is \textit{non-expansive} due to proposition \ref{combined-nonexpansiveness}, and is \textit{asymptotically regular} using corollary \ref{corollaryasymptoticregularity}.
\end{itemize}


\begin{figure}
    \centering
    \includegraphics[width=0.6\linewidth]{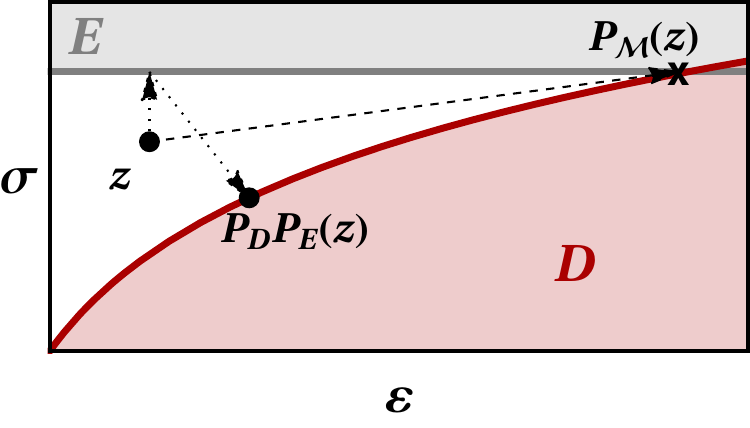}
    \caption{Simple illustration of the PSIs method: scheme of a one-element phase-space featuring the two subsets, see also the two projections that intervene in \cref{eq:convrate} to define the convergence rate. In this case, the convex sets are understood to be $D = \{ (\varepsilon', \sigma') \in Z : \sigma' - m(\varepsilon') \le 0 \}$ and $E = \{ (\varepsilon', \sigma') \in Z : \sigma' - F/A \ge 0 \}$. This setting was also considered in Ref.\,\cite{Ladeveze:2012}. 
}
    \label{fig:convex_sets}
\end{figure}

\paragraph{Step 3: apply Opial finite-dimensional theorem}
As the last step, having proved that the composition of projections $P(\boldsymbol{x})$ is a non-expansive, asymptotically regular map, it is possible to apply the finite-dimensional version of Opial's theorem (\ref{Opial}) to $P(\boldsymbol{x})$. 
Then, from proposition \ref{fixedpointsandsolution} one concludes that the algorithm converges to a point at the intersection between the physically-admissible set and the materially admissible set. 
\end{proof}



\subsection{Uniqueness of the solution}

The previous proof guarantees that the algorithm will converge to one point belonging to the intersection between the two closed convex subsets. 
This does not preclude that there might be infinite points belonging to this intersection, i.e., those where $E$ and $D$ overlap, recall \Cref{fig:convex_sets}. 
%
%
Using a zero-strain-zero-stress initialization we choose an initial phase-space point that for sure belongs to the frontier of the materially-admissible set and out of the physically-admissible one, thus the solution is expected to be unique and to be at the frontier (see \Cref{sec:influence_distance}, \Cref{fig:evolution_projections} in particular). 
%
%
Let us assume  henceforth that $\boldsymbol{z}^* = E \cap D$ is unique.

\subsection{Convergence of the algorithm}

%
We fetch MAP literature results \cite{MAP,lewismalick} for the case of two subspaces, and we will outline in the following how they apply in our case. 

%
Consider $\mathcal{H}$ to be a real Hilbert space, 
for a closed linear subspace $S$ of $\mathcal{H}$ we denote by $S^\perp$ its orthogonal complement in $\mathcal{H}$, and by $P_S$ the orthogonal projection of $\mathcal{H}$ onto $S$.
It was proved by J. von Neumann \cite{Neumann:1949} that for two closed subspaces $\mathcal{M}_1$ and $\mathcal{M}_2$ of $\mathcal{H}$, with intersection $M = \mathcal{M}_1 \cap \mathcal{M}_2$, the following convergence result holds:

\begin{align}
    \lim_{n \to \infty} \| (P_{\mathcal{M}_2} P_{\mathcal{M}_1})^n(\boldsymbol{x}) - P_\mathcal{M}(\boldsymbol{x}) \| = 0 \quad (x \in \mathcal{H}) \, ,
\end{align}

what is consistent with asymptotic regularity. 
This result says that the iterates $P_{\mathcal{M}_2} P_{\mathcal{M}_1}$ are strongly convergent to $P_\mathcal{M}$ \cite{convrate}.

Considering $2$ subspaces (i.e., the materially admissible set id the material was linear and an adjusted version of the physically-admissible set), a considerably detailed description of the rate of convergence is known, in terms of the notion of ``angle between subspaces'', a.k.a. the Friedrichs angle \cite{convrate}. 

\vspace{0.25cm}

\begin{definition}[Friedrichs angle $\theta$]
 Let $\mathcal{M}_1$ and $\mathcal{M}_2$ be two closed subspaces of the Hilbert space $\mathcal{H}$ with intersection $\mathcal{M} = \mathcal{M}_1 \cap \mathcal{M}_2$. The Friedrichs angle between the subspaces $\mathcal{M}_1$ and $\mathcal{M}_2$ is defined to be the angle $\theta \in [0, \pi/2]$ whose cosine is given by
\begin{equation}
    \cos \theta(\mathcal{M}_1, \mathcal{M}_2) 
    = 
    \sup\left\{ 
    \left| \langle \boldsymbol{x}, \boldsymbol{y} \rangle \right| 
    : 
    \boldsymbol{x} \in \mathcal{M}_1 \cap M^\perp \cap B_H, 
    \boldsymbol{y} \in \mathcal{M}_2 \cap M^\perp \cap B_H
    \right\}
    \label{friedrichsangle},
\end{equation}
where:
\begin{equation}
B_H = \{h \in \mathcal{H} : \|h\| \leq 1\}\label{ball}
\end{equation}
is the unit ball of $\mathcal{H}$.   
\end{definition}

It was proved as upper bound, and later as equality \cite{convrate}, that:
\begin{equation}
    \| (P_{\mathcal{M}_2} P_{\mathcal{M}_1})^n - P_\mathcal{M} \| 
    = 
    \left(
    \cos \theta(\mathcal{M}_1, \mathcal{M}_2)
    \right)^{2n-1} \, ,
    \quad 
    n \geq 1 \, .
    \label{eq:convrate}
\end{equation} 
This formula shows that the sequence  of iterates \( (P_{\mathcal{M}_2} P_{\mathcal{M}_1})^n \) converges uniformly to $P_M$  if and only if $\cos \theta < 1$ , i.e., if the Friedrichs angle between $\mathcal{M}_1$ and $\mathcal{M}_2$ is positive \cite{convrate}. 
When this happens, the iterates of $P_{\mathcal{M}_2} P_{\mathcal{M}_1}$ converge at geometrical rate to $P_\mathcal{M}$, in the following sense of ``quick uniform convergence'' \cite{convrate}:

\vspace{0.25cm}

\begin{proposition}[quick uniform convergence]
    $\text{There exist } A > 0 \text{ and } \beta \in (0, 1) \text{ such that }
    \| (P_{\mathcal{M}_2} P_{\mathcal{M}_1})^n - P_\mathcal{M} \| \leq A \beta^n \quad (n \geq 1).$
    \label{prop:quc}
\end{proposition}

Note that the MAP results presented in this subsection are formulated under the assumption that we are working with subspaces, even though the boundaries of neither $E$ nor $D$ directly satisfy this condition. Specifically, the former's is an affine subspace (i.e., it does not contain the origin), while the latter's, although it includes the origin (representing materials without residual stresses), is not closed under addition, recall \Cref{fig:convex_sets}. 
The geometric convergence rate is also valid in our case because, on the one hand, the affine subspace $E$ can always be converted into a subspace by translating the frame of reference, and, on the other hand, the convergence of the iterations over $D$ can be bounded by the convergence of iterations over a certain subspace., see \Cref{sec:bound}.
%

%


%% file: subfiles/1D_bar.tex
\subsection{Solving in closed form}

We resort to the closed-form expressions (that can be readily obtained in this simple setting) to clarify the concepts that we have been handling so far. 
Since there is only one element, the local and global phase spaces coincide, thus the subscript $e$ becomes redundant and will be dropped henceforth. 
This example is illustrative as it shows how the abstract concepts discussed in the previous section emerge naturally. 

As mentioned earlier, the starting point in phase space is $\boldsymbol{z'}_0 
= 
[ \varepsilon'_0 \quad \sigma'_0]^{\top} 
\in D$. 
The equilibrium equation, $\sigma A - F = 0 $, is indifferent to displacements. 
Compatibility says $\varepsilon = u / L$. 
The projection functional (\ref{eq:functional_Pi_E_general}) is simply 
\begin{align}
    \Pi_E
    =
    d^2 
    -
    \eta
    \left( 
    \sigma A - F
    \right)
    =
    w
    \left(
    { (\sigma - \sigma')^2 
    \over 
    2 C}
    +
    { (u/L - \varepsilon')^2 
    \over 
    2 C^{-1}}
    \right)
    -
    \eta
    \left( 
    \sigma A - F
    \right) \, ,
\end{align}
where $w = A L$ is the volume of the element. 
Enforcing stationarity of this functional, \cref{eq:stationarity}, leads to:

\begin{align}
    u = L \varepsilon'
    \Rightarrow 
    \varepsilon = \varepsilon'
    \, , \qquad
    \sigma
    =
    \sigma'
    + 
    C
    {\eta \over L}
    \, ,  \qquad
    \sigma 
    =
    F / A 
    \, .
\end{align}

Hence, $\boldsymbol{z}_0=P_E(\boldsymbol{z'}_0) \in E$ can be expressed as a linear operator $P_E : \mathbb{R}^2 \to \mathbb{R}^2$ such that 
\begin{align}
    \boldsymbol{z}_0
    =
    \begin{bmatrix}
    	\varepsilon_1 \\
        \sigma_1 
    \end{bmatrix}
    =
    P_E(\boldsymbol{z'}_0)
    =
    \begin{bmatrix}
    	1 &0 \\
        0 &0
    \end{bmatrix}
    \begin{bmatrix}
    	\varepsilon_0' \\
        \sigma_0' 
    \end{bmatrix}
    +
    \begin{bmatrix}
    	0 \\
        F/A 
    \end{bmatrix} \, .
\end{align}

%
%
From the projection onto $E$, we obtain stress-strain parts for each element, all of them satisfying equilibrium and compatibility. 
Next, project them back onto the closest point on their respective constitutive law: 
the projection functional (for a given $[\varepsilon \quad \sigma]^{\top}$) is given by \cref{eq:functional_Pi_D_general}, its E-L equations can be combined into \cref{eq:E-L_C_truss}, which, after using $d m / d \varepsilon = Y$ (linear-elastic constitutive lwaw), yields
\begin{align}
    \varepsilon' 
    =
    \varepsilon
    +
    {\left[ \sigma - Y \varepsilon' \right] 
    \over C^2}
    Y
    =
    \varepsilon
    +
    \left({ Y \over C} \right)^2
    {\sigma  \over Y}
    -
    \left({ Y \over C} \right)^2
    \varepsilon'
    \to
    \varepsilon'
    =
    { 1
    \over
    1 + \alpha^2}
    \varepsilon 
    +
    {\alpha^2 
    \over
    1 + \alpha^2}
    {\sigma \over Y}
    \, ,
    \label{eq:partial_1}
\end{align}
where $\alpha = Y / C$.
From the constitutive law, $\sigma' = Y \varepsilon'$,       
\begin{align}
    \sigma'
    =
    { 1
    \over
    1 + \alpha^2}
    Y \varepsilon 
    +
    {\alpha^2 
    \over
    1 + \alpha^2}
    \sigma
    \, ,
    \label{eq:partial_2}
\end{align}
Gathering \cref{eq:partial_1,eq:partial_2}, one can write
\begin{align}
    \boldsymbol{z}'_1
    &=
    P_D(\boldsymbol{z}_1)
    =
    \begin{bmatrix}
    	\varepsilon_1' \\
        \sigma_1' 
    \end{bmatrix}
    =
    {1 \over 1 + \alpha^2}
    \begin{bmatrix}
        1 &\alpha^2/Y \\
        Y &\alpha^2
    \end{bmatrix}
    \begin{bmatrix}
        \varepsilon_1 \\
        \sigma_1 
    \end{bmatrix}
    =
    {1 \over 1 + \alpha^2}
    \begin{bmatrix}
        1 &\alpha^2/Y \\
        Y &\alpha^2
    \end{bmatrix}
    \left(
    \begin{bmatrix}
    	1 &0 \\
        0 &0
    \end{bmatrix}
    \begin{bmatrix}
        \varepsilon_0' \\
        \sigma_0' 
    \end{bmatrix}
    +
    \begin{bmatrix}
    	0 \\
        F/A 
    \end{bmatrix}
    \right) \nonumber \\
    &=
    P_D P_E (\boldsymbol{z}'_0)
    =
    {1 \over 1 + \alpha^2}
    \begin{bmatrix}
        1 &0 \\
        Y &0
    \end{bmatrix}
    \begin{bmatrix}
        \varepsilon_0' \\
        \sigma_0' 
    \end{bmatrix}
    +
    {\alpha^2 \over 1 + \alpha^2}
    \begin{bmatrix}
        Y^{-1}F/A \\
        F/A 
    \end{bmatrix}
    \, .
\end{align}
This completes the first iteration. It thus follows that after $n$ iterations: 
\begin{align}
    \boldsymbol{z}'_n
    &=
    (P_D P_E)^n (\boldsymbol{z}'_0)
    =
    {1 \over (1 + \alpha^2)^n}
    \begin{bmatrix}
        1 &0 \\
        Y &0
    \end{bmatrix}
    \begin{bmatrix}
        \varepsilon_0' \\
        \sigma_0' 
    \end{bmatrix}
    +
    \alpha^2
    \left(
    \sum_{i=1}^{n}
    {1 \over (1 + \alpha^2)^i}
    \right)
    \begin{bmatrix}
        Y^{-1}F/A \\
        F/A 
    \end{bmatrix}
    \, .
    \label{eq:explicititeration}
\end{align}




We are about to derive explicitly the convergence rate in the case of two convex sets defined by linear equations in 1D for one element. 
%
%
The solution in this simple case is known beforehand: $\boldsymbol{z}^* = [Y^{-1}F/A \quad F/A]^\top$. 

We want to estimate at iteration $n$ the error given by $\boldsymbol{z}^*-\boldsymbol{z}'_n$. 
From \cref{eq:explicititeration}, choose as initial point for the algorithm the vector $[\varepsilon_0' \quad \sigma_0']^\top =
[ 0 \quad 0 ]^\top$ 
which is always a valid point belonging to the constitutive law but note also that the choice of starting points is irrelevant, what becomes obvious if in \cref{eq:explicititeration} one takes the limit $n \to \infty$.
Then:
\begin{align}
    \boldsymbol{z}^*-\boldsymbol{z}'_n
    =
    \left(1-\alpha^2
    \sum_{i=1}^{n}
    {1 \over (1 + \alpha^2)^i}
    \right)
    \begin{bmatrix}
        Y^{-1}F/A \\
        F/A 
    \end{bmatrix}= \left(1-\frac{\alpha^2}{1+\alpha^2}
    \sum_{i=1}^{n}
    {1 \over (1 + \alpha^2)^{i-1}}
    \right)
    \begin{bmatrix}
        Y^{-1}F/A \\
        F/A 
    \end{bmatrix}
    \, .
\end{align}

Now, with $q:=1/({1+\alpha^2})$, recall that for  partial geometric series:
$
    \sum_{i=0}^{n-1} q^i=\sum_{i=1}^{n} q^{i-1}=\frac{1-q^n}{1-q}  \, .
$
We use this property to derive:

\begin{align}
    \boldsymbol{z}^*-\boldsymbol{z}'_n
    =
    \left(
    1-\frac{\alpha^2}{1+\alpha^2}
    {{1-\frac{1}{(1 + \alpha^2)^{n}}}\over {1-\frac{1}{(1 + \alpha^2)}}}
    \right)
    \begin{bmatrix}
        Y^{-1}F/A \\
        F/A 
    \end{bmatrix}
    = 
    \frac{1}{(1 + \alpha^2)^{n}}
    \begin{bmatrix}
        Y^{-1}F/A \\
        F/A 
    \end{bmatrix}
    \,
\end{align}

and taking the phase space norm for both terms leads to:

\begin{align}
    \| \boldsymbol{z}^*-\boldsymbol{z}'_n  \| 
    &= 
    \left(
    \frac{1}{1+\alpha^2}
    \right)^n 
    \| \boldsymbol{z}^* \| 
    = 
    (\cos \theta)^{2n} 
    \| \boldsymbol{z}^* \| \, ,
    \label{practical}
\end{align}
after setting $\alpha=\tan \theta$ and using the trigonometric identity $(1+{(\tan \theta)}^2)^{-1}={(\cos \theta)}^2$. 
This shows how the angle $\theta$ (the Friedrichs angle), which depends on $\alpha$ and on $C$ as a consequence, affects the convergence rate, 
consistently with \cref{eq:convrate}. 
Revisit also \Cref{fig:iterations_PS,fig:iterations_PS_normalized}.


From this explicit computation, 
we expect the constant which defines the convergence rate to be equal to $(1+\alpha^2)^{-1}$. 
For illustration purposes, let us recover this result directly from the definition \ref{friedrichsangle}. 
The scalar product 
is the conventional one in $\mathbb{R}^2$
, and
%
\[
\mathcal{M}_1:=\left\{\sigma - F/A = 0, \forall \varepsilon \right\}, \qquad
\mathcal{M}_2:=\left\{\sigma - Y\varepsilon = 0, \forall \varepsilon \right\}
\]
but, in this form, $\mathcal{M}_1$ is not a subspace. 
However, we can translate the horizontal axis, what does not change the definition of angle, and rewrite:
\[
\mathcal{M}_1:=\left\{\Delta \sigma = \sigma - F/A = 0, \forall \varepsilon \right\},
\]
which does represent a subspace (namely, the horizontal axis itself). 
%
%
Now, 
$M = \mathcal{M}_1 \cap \mathcal{M}_2$=$\left\{{0} \right\}$,
hence $M^\perp=\left\{v: v \neq {0} \right\}$. 
Defining vectors  
$\boldsymbol{x} \in \mathcal{M}_1$ 
such that 
$\left\| \boldsymbol{x} \right\|=1 $ 
and
$\boldsymbol{y} \in \mathcal{M}_2$ such that $ \left\| \boldsymbol{y} \right\|=1$, 
we end up with:
\begin{align*}
    \boldsymbol{x}
    =\begin{pmatrix}
    \sqrt{\frac{2}{C w_e}} \\
    0
    \end{pmatrix} \qquad
    \boldsymbol{y} 
    =\begin{pmatrix}
    \sqrt{\frac{2C}{(C^2+Y^2) w_e}} \\
    Y\sqrt{\frac{2C}{(C^2+Y^2) w_e}}
    \end{pmatrix} \, .
\end{align*}

The Friedrichs angle is computed using the scalar product as 
$
 \left| \langle \boldsymbol{x}, 
 \boldsymbol{y} \rangle \right|
 =(1+\alpha^2)^{-1/2}\,
$.
Finally, exploiting \cref{eq:convrate}, we have that the convergence rate is defined by the square of the previous quantity, namely by $(1+\alpha^2)^{-1}$, the same result as derived in explicit computation. 
%

\begin{figure}
    \centering
    \includegraphics[width=0.6\linewidth]{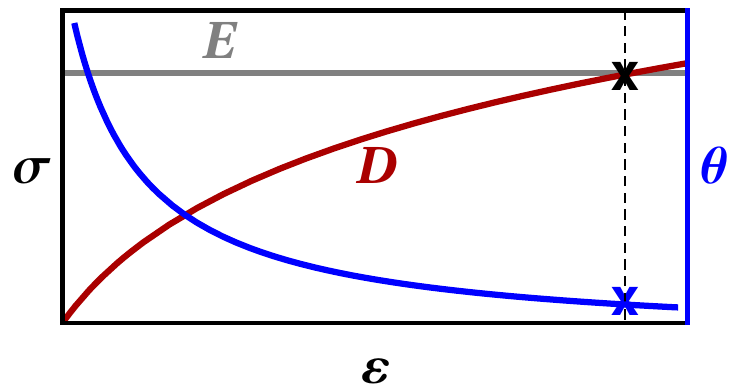}
    \caption{Sketch of the Friedrichs angle in the 1D one-element case. The left vertical axis depicts stresses, while the right vertical one depicts the slope angle of the constitutive; both share horizontal axis (strains). The black cross marks the solution, while the blue one marks the ``local'' Friedrichs angle value controlling the performance of the PSI method in this example.}
    \label{fig:angles_plot}
\end{figure}

\subsection{Adding non-linearity complicates the analysis}
\label{sec:non-linearity}

Let us now turn to considering the effects of material non-linearity. 
When a generic non-linear constitutive law is considered, it can be rewritten through Taylor expansion up to the second order to have a better approximation with respect to a linear approximation. 
Infinitesimal strains $\varepsilon$ are still in use. 
In the 1D case, the non-linear constitutive law is approximated around to $\varepsilon = 0$ as
$\sigma(\varepsilon) \approx a + b \varepsilon + c \varepsilon^2$,
where $a=\sigma(0), b=\left. \frac{{d\sigma}}{{d\varepsilon}} \right|_{\varepsilon=0},  c=  \frac{{1}}{{2}} \left. \frac{{d^2\sigma}}{{d\varepsilon^2}} \right|_{\varepsilon=0} $. 
Under the assumption $a=\sigma(0)=0$ (no residual stresses) and considering that the equation must define the convex set, the signs in the equation would be $b>0, c<0$. 
Since the case considered is the one with a small non-linear term, we can rewrite the equation as:

\begin{equation}
\sigma(\varepsilon) = Y(\varepsilon - k\varepsilon^2) \, ,
\label{eq:constitutive_law_NL_k}
\end{equation}

assuming $k\ll1$, so that in the limit as $k \to 0$ a linear constitutive law with Young modulus $Y$ is recovered. 
To study the convergence behavior, as before, we will focus on the case of given level load.
The projection onto the physically-admissible set is such that
$\boldsymbol{z}_{n+1}=P_E(\boldsymbol{z'}_n) \in E$ can be expressed as a linear operator $P_E : \mathbb{R}^2 \to \mathbb{R}^2$ as:
\begin{align}
    \boldsymbol{z}_{n+1}
    =
    \begin{bmatrix}
    	\varepsilon_{n+1} \\
        \sigma_{n+1}
    \end{bmatrix}
    =
    P_E(\boldsymbol{z'}_n)
    =
    \begin{bmatrix}
    	1 &0 \\
        0 &0
    \end{bmatrix}
    \begin{bmatrix}
    	\varepsilon_n' \\
        \sigma_n' 
    \end{bmatrix}
    +
    \begin{bmatrix}
    	0 \\
        F/A 
    \end{bmatrix} \, .
\end{align}

while the projection onto the materially admissible set can be performed through \cref{eq:E-L_C_truss},
which for the non-linear constitutive law \cref{eq:constitutive_law_NL_k},
after some simplifications, 
boils down to
\begin{equation}
    \varepsilon'^{3}(-2\alpha^2k^2)
    +
    \varepsilon'^2(3k\alpha^2)
    +
    \varepsilon'(-1-\frac{2kY\sigma}{C^2}-\alpha^2)
    +
    (\varepsilon+\frac{\sigma Y}{C^2})
    =
    0
\label{eq:3rdordereq}
\end{equation}

where $k, \varepsilon, \sigma, \alpha, C, Y$ are all known quantities, $\varepsilon'$ represents the solution to be found. 
Since both the term $\varepsilon'^3 $
and $\varepsilon'^2$ disappear for $k=0$, we are looking for roots of an algebraic equation under a \textit{singular perturbation} \cite{perturbationbook}. 
From the fundamental theorem of algebra, we only have two cases for a real third order algebraic equation: 
(i) three real roots, 
(ii) one real roots and two imaginary roots. 
Since we are not guaranteed to have a unique real root converging to the solution when $k=0$, we would need to study qualitatively the behavior of such roots. 
An approximate approach is looking for a \textit{regular expansion} of the solution in terms of the small parameter:

\begin{equation}
    \varepsilon'
    =
    \varepsilon_{0}'
    +
    k\varepsilon_{1}'
    +
    \mathcal{O}(k^2) \, 
\label{kexpansion}
\end{equation}

We can substitute into eq. (\ref{eq:3rdordereq}) and apply the fundamental theorem of perturbation theory \cite{perturbationbook}.
The following approximated expression can be derived:
\begin{align}
     \varepsilon_{n}'
    =
    \left(
    { 1
    \over
    1 + \alpha^2}
    \varepsilon 
    +
    {\alpha^2 
    \over
    1 + \alpha^2}
    {\sigma \over Y}
    \right)
    +k
    \left( 
    \frac{Y(C^2\varepsilon+\sigma Y)(\sigma Y^2 +C^2(-2 \sigma +3 \varepsilon Y))}{(C^2+Y^2)^3}
    \right) \, .
\end{align}
The extra term renders the new expression more difficult to handle, thus a clean recurrence relation could not be found. 
This led us to consider an alternative path: finding a result for the non-linear case in terms of the one for the linear, \cref{eq:explicititeration}.





\subsection{A general bound for the convergence rate based on linear analysis}
\label{sec:bound}

\begin{figure}
\centering
\captionsetup[subfigure]{justification=centering}
\begin{subfigure}[t]{.45\linewidth}
\includegraphics[width=\linewidth]{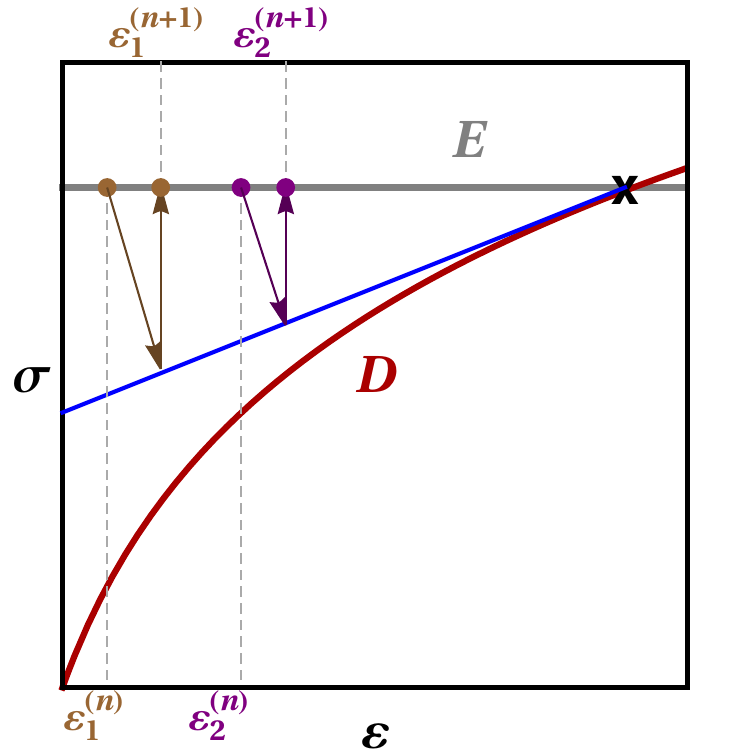}
\caption{\Large (a)}\label{fig:observation_i}
\end{subfigure}
\begin{subfigure}[t]{.45\linewidth}
    \includegraphics[width=\textwidth]{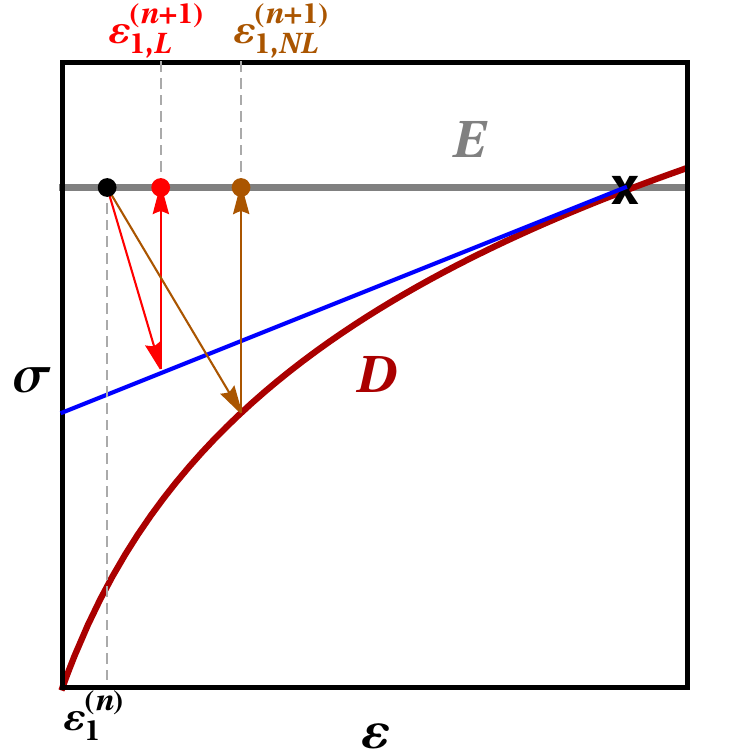}
    \caption{\Large (b)}
    \label{fig:observation_ii}
\end{subfigure}
        \caption{Visual representations of the key observations. 
        (a) First: projections onto linear envelope preserve ordering: if $\varepsilon_1^{(n)}<\varepsilon_2^{(n)}$ then $\varepsilon_1^{(n+1)}<\varepsilon_2^{(n+1)}$. 
        (b) Second: projection onto the actual constitutive law move faster to the solution than the projection onto the linear envelope: $\varepsilon_{1,L}^{(n+1)}<\varepsilon_{1,NL}^{(n+1)}$. 
        }
\end{figure}

%
We aim for an upper bound for the ``error'' of the algorithm when applied to a non-linear case, in order to estimate the convergence rate based on the result for the linear one. 
The solution is given by $
\boldsymbol{z}^*=[\varepsilon^* \quad \sigma^*]^{\top}$. 
Inspired by Ref. \cite{lewismalick}, 
the key is to observe that since the constitutive law is convex, at $\boldsymbol{z}^*$ the slope (angle) of the constitutive law is greater than the slope of the tangent at $\boldsymbol{z}^*$ for any $\varepsilon<\varepsilon^*$, regard \Cref{fig:angles_plot}. 
%
Recall -- \cref{eq:convrate} -- that \textit{the smaller the Friedrichs angle, the slower the convergence}. 
It is possible to show that this \textit{smallest} value of the angle between sets $D$ and $E$, given by the angle computed at $\boldsymbol{z}^*$, represents a bound for the actual convergence of the non-linear iterations. 
We define the ``bounding line'' to be the tangent at $\boldsymbol{z}^*$. 
%
%
This is explained in \Cref{fig:observation_i,fig:observation_ii}. 
In the limit case: the constitutive law would feature a maximum and the equilibrium set would be tangent to it at a point; as shown with an analogous example with a parabola in Ref. \cite{gubin1966method}, the convergence rate is known to be slower than any geometric progression in such a scenario. 
%

Two fundamental observations are necessary in order to formalize this analysis. 
In the following, we will consider $P(\boldsymbol{z})=P_{E}(P_{D}(\boldsymbol{z}))$ and $P_L(x)=P_{E}(P_{D,L}(\boldsymbol{z}))$, where $P_{D,L}$ is a similar projection but onto the bounding line, and the subscripts are used not to denote different elements but different illustrative points. 

\paragraph{Observation (i)} See \Cref{fig:observation_i}. 
For a convex constitutive law, 
given two points $\boldsymbol{z}^{(n)}_1 = [\varepsilon_1 \quad \sigma_1]^{\top} \, , \boldsymbol{z}^{(n)}_2 = [\varepsilon_2 \quad \sigma_1]^{\top} \in E$, 
if $\varepsilon^{(n)}_1<\varepsilon^{(n)}_2$, then 
defining $\boldsymbol{z}^{(n+1)}_{1} = P_L(\boldsymbol{z}^{(n)}_1)$ 
and 
$\boldsymbol{z}^{(n)}_{2} = P_L(\boldsymbol{z}^{(n)}_2)$, 
it follows also that 
$\varepsilon^{(n+1)}_{1}<\varepsilon^{(n+1)}_{2}$. 
In other words, the ordering is preserved. 
 
%



\paragraph{Observation (ii)} 
See \Cref{fig:observation_ii}. 
Given a point $\boldsymbol{z}^{(n)}_1 = [\varepsilon^{(n)}_1 \quad \sigma^{(n)}_1]^{\top} \in E$, 
taking $\boldsymbol{z}^{(n+1)}_1$ as $P(\boldsymbol{z}^{(n)}_1)$, 
which is a standard projection onto the constitutive law, 
and $\boldsymbol{z}^{(n+1)}_{1,L}$ which is a projection onto the bounding line, 
then $\varepsilon^{(n+1)}_{1,L}<\varepsilon^{(n+1)}_{1}$. 

\paragraph{From observations to an upper bound for the error}

Combining observations (i) and (ii), it follows that $ \| \boldsymbol{z}^* - \boldsymbol{z}^{(n+1)} \| < \| \boldsymbol{z}^* - \boldsymbol{z}_L^{(n+1)} \| $, but now we recognize that we can apply \Cref{prop:quc} to the right-hand side of the inequality, thus reaching finally that $\| \boldsymbol{z}^* - \boldsymbol{z}^{(n+1)} \| = \mathcal{O}(\beta^n)$, for some $\beta<1$. 
In summary, this demonstrates that the convergence rate of the PSI method is geometric as long as the constitutive law and the physical-admissibility conditions (equilibrium and compatibility) render two convex sets. 

This concludes the analysis of this example and the discussion as to convergence properties. 
The next section tackles a more realistic engineering example which serves to illustrate the method in an application. 

\subsection{Future steps: generalizing}

The previous section was limited to the simplest system (one element) and the simplest type of element (uniaxial stress bar). 
For this case, $N_e=n_e=1$, the concepts borrowed from MAP and POCS are applicable with no tweak whatsoever. 
The projection operator \cref{eq:P_E} does not meet automatically the definition \ref{thm:HilbertProjection} in a larger global phase-space $Z$, since it is unclear under what conditions it can be considered as the frontier of a ``closed convex set''. 
By way of illustration, consider two bars ($n_e = 1$ but now $N_e = 2$), hence \cref{eq:equilibrium} represents a space of co-dimension 2 in $\mathbb{R}^4$: Can it be the frontier of a closed convex set? 
Can an auxiliary closed convex set be defined naturally at least? 
Answering these questions would clarify the geometric picture and hopefully would unravel how to generalize the use of alternating projections (to enable proving existence) and the of Friedrichs angle (to characterize the convergence rate).  

%% file: subfiles/Kirchdoerfer.tex

\begin{figure}
\centering
\captionsetup[subfigure]{justification=centering}
\begin{subfigure}[t]{.45\linewidth}
\includegraphics[width=\linewidth]{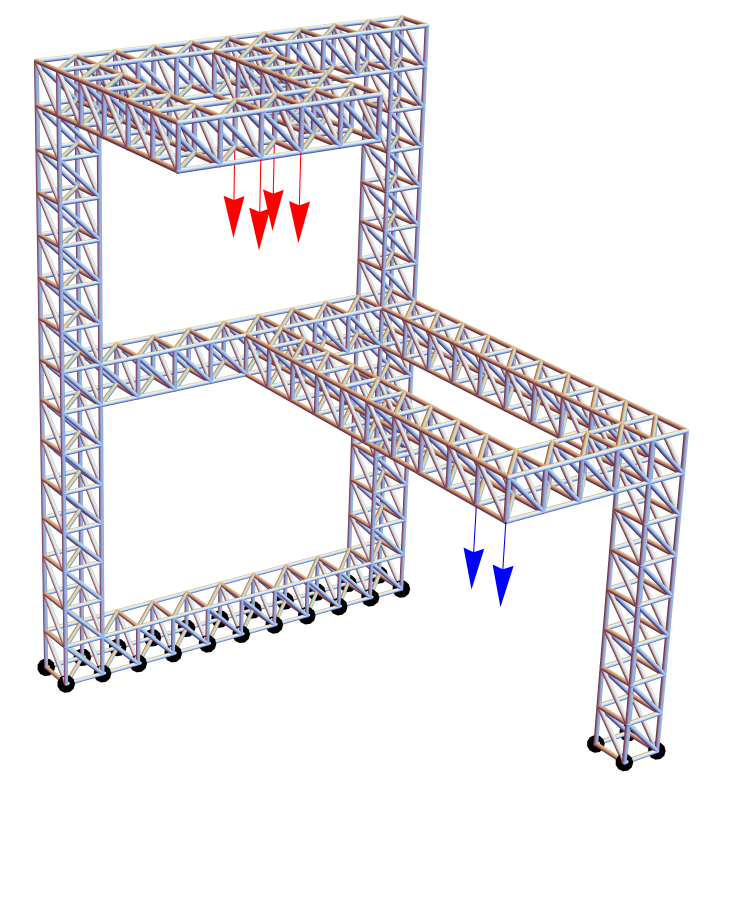}
\caption{\Large (a)}\label{fig:truss_bars}
\end{subfigure}
\begin{subfigure}[t]{.54\linewidth}
\includegraphics[width=0.99\textwidth]{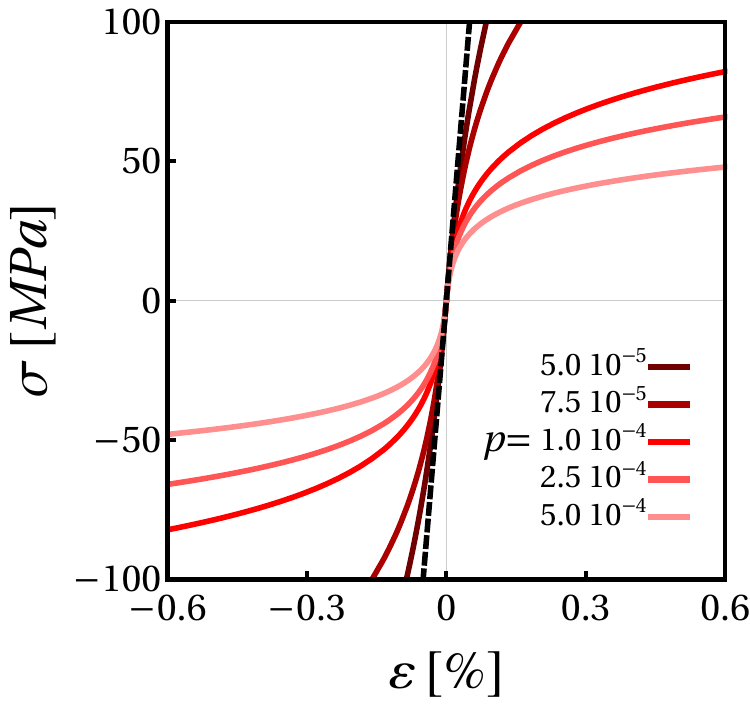}
    \caption{\Large (b)}
    \label{fig:1D_constitutive_law}
\end{subfigure}
        \caption{Truss example. (a) Structure: red arrows indicate applied forces, blue ones correspond to imposed displacements. (b) Constitutive law used for the all bars in the truss.}
\end{figure}

The structure we focus on is ``Kirchdoerfer's truss'', a simple truss system introduced in Ref.\,\cite{Trent_1}. 
It features 1246 bars (elements) and 376 connections (nodes), see \Cref{fig:truss_bars}.
The non-linear elastic law presented in \cite{Trent_1} is equipped on all the elements, and it is given by
\begin{equation}
\label{eq:material_law}
    \sigma
    =
    m(\varepsilon)
    =
    Y_0 
    \left[
    (|\varepsilon|+c)^p - c^p
    \right]
    \mathrm{sign}(\varepsilon) \, ,
\end{equation}
where $c = p^{(1-p)^{-1}}$, $Y_0$ represents the axial stiffness at zero strain, $|\cdot|$ is the absolute value function and $\mathrm{sign}(\cdot)$ is the sign function, 
see \Cref{fig:1D_constitutive_law};, 
in this example, $Y_0 = 200 \mathrm{GPa}$ and $p$  will be varied. 
There are four nodes where external vertical forces are applied (magnitudes -1000, -1000, -100 and +1800 newtons, see supplementary material for complete details), and two nodes where a vertical downward 5-cm displacement is enforced as boundary condition.   

We compare the new solver to a damped Newton-Raphson (see appendix), which is more adequate to handle the forced displacement conditions. 
We will use a damping parameter of $ 0.8$, i.e., the zero-strain stiffness is 20\% of each iteration stiffness matrix. 

The force residual tolerance for NR and PSI is initially set to be $\mathrm{tol}_1 =5 \cdot 10^{-2} = 5 \% $, the phase-space distance tolerance (recall \Cref{sec:consecutive_iterations}) is $\mathrm{tol}_2 = \mathrm{tol}_1/10$. 

Both the PSI and NR methods were implemented by the same user (the last author) in Mathematica notebooks \cite{Mathematica} (version 12.2), made available via GitHub (see supplementary material). 
Neither implementation has been optimized by a computer science expert; 
a comprehensive comparison on a dedicated solid mechanics open-source software (e.g., Akantu \cite{Akantu}, Deal II \cite{Deal:II} or FEniCS \cite{FENICS:book}) is out of the scope of this paper that focuses on presenting the method and justifying its mathematical properties. 

All the execution times correspond to the average of five runs. 
All the simulations were run in the same machine with specifications: 
\begin{itemize}[leftmargin=*]
    \item CPU: 11th Gen Intel(R) Core(TM) i7-1165G7 @ 2.80GHz, 4 cores (Mathematica runs one ``kernel'' -- task -- per core, using a single thread per core). 
    \item RAM: Intel(R) Tiger Lake-LP Shared SRAM (rev 20), 15Gi total. 
    \item OS: Ubuntu 22.04.4 LTS.
\end{itemize}


\subsection{Detailed analysis of one run}

The material model corresponds to $p=10^{-4}$, see \Cref{fig:1D_constitutive_law}. 
We initially employ phase-space distance minimization for projections onto $D$, i.e., 
we are using Mathematica's function \texttt{FindMinimum} to solve \cref{eq:functional_Pi_D_general}, 
with the ``principal axis'' (PA) method, 
step control (``trust region'' type), 
while also providing the previous physically-admissible strain ($\boldsymbol{\varepsilon}_e$) as starting point to the search.

\begin{figure}
\centering
\includegraphics[width=0.85\textwidth]{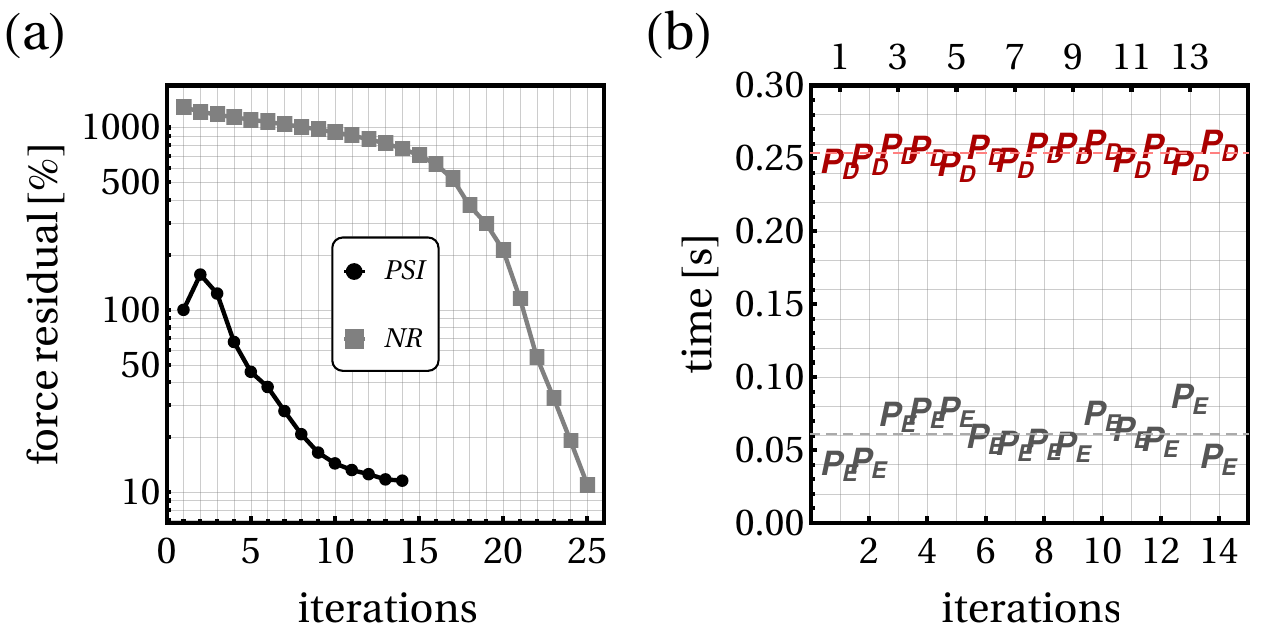}
        \caption{Truss example: (a) Evolution of the residual over iterations, Newton-Raphson solver (NR) v. phase-space iterative (PSI) solver  (b) PSI time breakdown in each iteration: time invested in performing projection onto physically-admissible set ($P_E$) and onto materially-admissible set ($P_D$).}
        \label{fig:NR_PSI_run_comaprison}
\end{figure}

\Cref{fig:NR_PSI_run_comaprison}(a) displays two interesting traits of the PSI when compared to NR. 
On one hand, the force residual can increase after some iterations (such is the case after the second one). 
On the other hand, the convergence can be very fast in NR, once it closes in the solution, but phase-space distance can reach a certain level of precision in fewer iterations. 

Let us analyze the time breakdown of one PSI run, \Cref{fig:NR_PSI_run_comaprison}(b). 
There are two tasks per iteration: 
solving two linear systems (\cref{eq:sol_eta,eq:sol_u} for $\boldsymbol{u}$ and $\boldsymbol{\eta}$) to then evaluate $\boldsymbol{\sigma}$ and $\boldsymbol{\varepsilon}$ in the projection onto $E$ (\cref{eq:compatibility,eq:sigma_update_P_E} operatively), 
followed by getting $\boldsymbol{\sigma}'$ and $\boldsymbol{\varepsilon}'$, 
\cref{eq:compatibility}, 
through the projection onto $D$, $P_D$, which is performed by means of local phase-space distance minimization, 
i.e., solving either a non-linear equation for each element \cref{eq:E-L_C_truss} (note that in this case it is a scalar equation, but will not be the case in general) or minimizing a distance.  
\Cref{fig:NR_PSI_run_comaprison}(b) consistently shows that the $P_D$ step consumes about five times as much time as $P_E$. 
In absolute terms, $P_D$ takes $0.25$ seconds on average, while $P_E$ takes $0.06$.  
Logically, the former task would be the first to be dealt with if it came to reducing time, 
what could be done in two ways: 
(a) by improving on the performance of the distance-minimization method, 
(b) by improving on the parallelization, i.e., bringing more processors to share in the burden of solving one non-linear algebraic equation per element. 

\subsection{Influence of distance choice}
\label{sec:influence_distance}

%
The results using the implementation of the NR method are used as ground truth to compare. 

\begin{figure}
    \centering    
    \includegraphics[width=\textwidth]{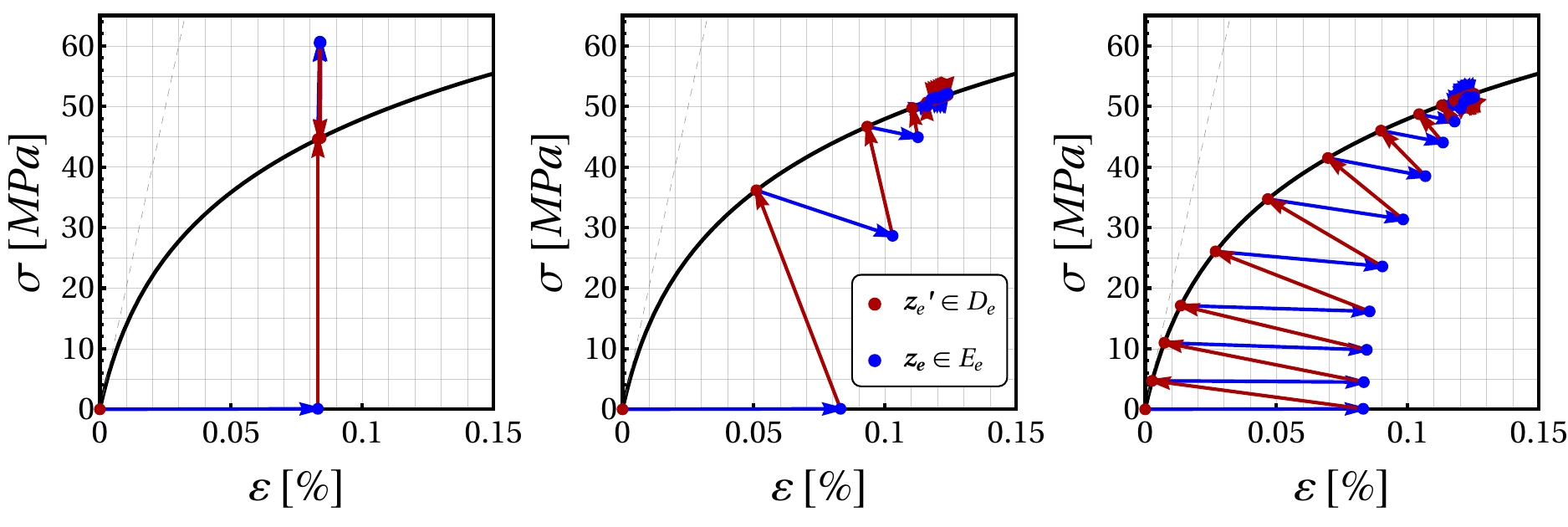}
    \includegraphics[width=\textwidth]{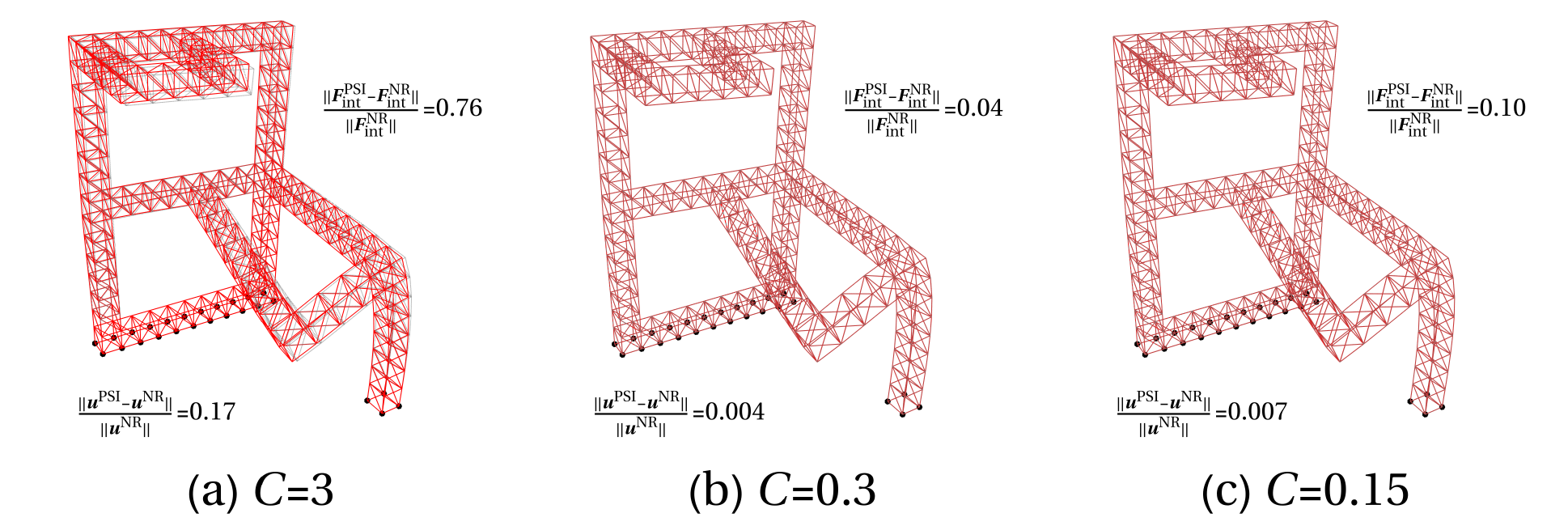}
    \caption{Iterations for one element, arrows represent successive projections, using distance minimization with $C$ equal to (a) $3.0 Y_0$ (left column), (b) $0.3 Y_0$ (center) and (c) $0.15 Y_0$ (right). 
    First projection yields strain to satisfy essential BCs. 
    The position of the physically-admissible set evolves as $E_e$ depends on $[\sigma_e]_{e=1}^{N_e}$, which change each time $P_E$ is applied at the structure level.}
    \label{fig:evolution_projections}
\end{figure}

\Cref{tab:table_1} along with \Cref{fig:evolution_projections} shows the effect of $C$, which was also sketched in \Cref{fig:first}. 
We remark, first, that all the simulations converge because the phase-space distance between consecutive points stops changing, recall the dual condition discussed in \Cref{sec:consecutive_iterations}. 
Choosing $C = 3 Y_0$ yields misleading fast convergence (only 1.6 seconds, taking 5 iterations): the phase-space distance criterion ceases to improve but the physically-admissible states are far from being reconciled with the constitutive law, as it is clearly seen in \Cref{sec:consecutive_iterations}(a) (top panel). 
Both the displacements and the internal forces are greatly misrepresented, the difference with respect to NR being $17 \%$ and $76\%$ respectively; the bottom panels of \Cref{sec:consecutive_iterations} displays the deformed shapes in each case.

\bgroup
\def\arraystretch{1.5}
\begin{table}[H]
\centering
\caption{Wallclock time comparisons using different values of $C$ and distance minimization (principal axis method) for $P_D$. Material defined by $p = 1.0 \cdot 10^{-4}$.}
\begin{tabularx}{0.95\textwidth}{ c| *{5}{Y|} }
\cline{2-5}
&\multicolumn{4}{|c|}{$C/Y_0$} \\
\cline{2-5}
& $3$ & $1$ &$0.3$  &$0.15$  \\ 
\hline
\multicolumn{1}{|c|}{Running time [s] \textbf{/} \# iterations}& \multicolumn{1}{c|}{1.6\textbf{/}5}  & 4.6\textbf{/}15   & 4.6\textbf{/}14  &10.6\textbf{/}35 \\ \hline
\end{tabularx}
 \label{tab:table_1}
\end{table}
\egroup

Reducing the value to $C = 1 Y_0$ increases the number of iterations with no improvement in the performance. 
Conversely, $C = 0.3 Y_0$ takes a similar amount of both times and iterations but radically improves the method performance: 
the bottom panel of \Cref{sec:consecutive_iterations}(b) reveals two almost-indistinguishable deformed shapes while the top panel explains why the method is so efficient and does not get stuck in this case. 
The fact that by further reducing $C = 0.15 Y_0$ penalizes the performance (both in terms of running time/number of iterations and of resemblance to the actual --NR-- solution) confirms the existence of an optimal distance function for operational purposes. 

\subsection{Influence of method to perform projection onto D}

Interestingly, cf. \Cref{tab:table_7}, the derivative-less method for solving this scalar minimization problem is faster than those that compute gradients. 
All methods take 14 iterations to converge, the final force imbalance is $6 \%$ in all cases (to compare with the NR convergence tolerance, set to be $5 \%$).

\bgroup
\def\arraystretch{1.5}
\begin{table}[H]
\centering
\caption{Wallclock time comparisons using distance minimization for $P_D$ (different methods). Material defined by $p = 1.0 \cdot 10^{-4}$. Distance defined by $C=0.3$. }
\begin{tabularx}{0.95\textwidth}{ c| *{5}{Y|} }
\cline{2-5}
& Principal Axis (derivative-free) & Conjugate Gradient &Quasi-Newton (BFGS)  &Newton  \\ \hline
 \multicolumn{1}{|c|}{Running time [s]}& \multicolumn{1}{c|}{4.6}  & 6.9   & 6.9  &6.3 \\ \hline
\end{tabularx}
 \label{tab:table_7}
\end{table}
\egroup

BFGS refers to Broyden–Fletcher–Goldfarb–Shanno method, a quasi-Newton method that does not require re-computing the Hessian \cite{FBGS:1989}. 


One can also resort to the corresponding E-L equations, \cref{eq:E-L_C_truss}, instead of going for the direct minimization when it comes to performing projections onto $D$. 
\Cref{tab:table_E-L} suggests that, at least in this case, solving directly the non-linear equation derived from the optimality condition is more time-consuming than performing direct minimization of the objective function.  

\bgroup
\def\arraystretch{1.5}
\begin{table}[H]
\centering
\caption{Wallclock time comparisons using distance minimization for $P_D$ (different methods). Material defined by $p = 1.0 \cdot 10^{-4}$. Distance defined by $C=0.3$. }
\begin{tabularx}{0.95\textwidth}{ c| *{4}{Y|} }
\cline{2-4}
& Newton &ACN  &Secant  \\ 
\hline
 \multicolumn{1}{|c|}{Running time [s]}& \multicolumn{1}{c|}{7.8}  & 28.7   & 9.2   \\ \hline
\end{tabularx}
 \label{tab:table_E-L}
\end{table}
\egroup

``ACN'' refers to Affine Covariant Newton method 
\cite{Deuflhard:2011}.

\subsection{Influence of number of kenerls}

\bgroup
\def\arraystretch{1.5}
\begin{table}[H]
\label{tab:table_8}
\centering
\caption{Wallclock time comparisons using distance minimization for $P_D$. Material defined by $p=1 \cdot 10^{-4}$. Using $C = 0.3$. Minimization approach with PA method.}
\begin{tabularx}{0.95\textwidth}{ c| *{5}{Y|} }
\cline{2-5}
&\multicolumn{1}{c}{\multirow{2}{*}{Newton-Raphson solver}}& \multicolumn{3}{|c|}{Phase-space iterative solver} \\ \cline{3-5}
 & & \multicolumn{1}{c|}{serial} & \multicolumn{1}{c|}{2 kernels} & \multicolumn{1}{c|}{4 kernels}  \\ \hline
 \multicolumn{1}{|c}{Running time [s]}& \multicolumn{1}{|c|}{1.8}  & 4.6   & 2.5   & 2.2 \\ \hline
\end{tabularx}
\end{table}
\egroup

The scaling is good up for 2 processors; 
afterward, Mathematica does not take advantage of the extra resources. 
The fastest case (4 kernels) is still a 10\% slower than NR.

\subsection{How does non-linearity affect relative performance?}

We observe that PSI fares better than NR when most material response is non-linear, i.e., low value of $p$ (see \Cref{tab:table_p}). 
This is attributed to PSI enforcing the constitutive law during a dedicated step ($P_D$) and element-wise, while NR enforces equilibrium and the constitutive law simultaneously when solving the residual equations. 
NR relies on the tangent stiffness matrix, so there being elements that accumulate deformation can lead to entries of disparate magnitude and poor conditioning of the tangent-stiffness matrices. 

\bgroup
\def\arraystretch{1.5}
\begin{table}[H]
\centering
\caption{Wallclock time comparisons using distance minimization for $P_D$.  Minimization approach with PA method.}
\begin{tabularx}{0.75\textwidth}{ c| *{4}{Y|} }
\cline{2-4}
&\multicolumn{3}{c|}{Degree of non-linearity, $p/10^{-4}$} \\ 
\cline{2-4}
 &  \multicolumn{1}{c|}{2.0 (stiff)} & \multicolumn{1}{c|}{1.0} & \multicolumn{1}{c|}{0.5 (compliant)} \\ \hline
 \multicolumn{1}{|c|}{NR running time [s]}& 1.08  & 1.8 & 3.1 \\ \hline
 \multicolumn{1}{|c|}{PSI running time [s]}&  1.9 ($C/Y_0 = .4 $)  & 2.2 ($C/Y_0 = .3 $) & 2.6 ($C/Y_0 = .2 $) \\ \hline
\end{tabularx}
\label{tab:table_p}
\end{table}
\egroup


\Cref{tab:table_p} suggests that PSI can be especially competitive when it comes to finding approximate solutions of problems featuring widespread non-linear behavior. 
%
%
The method converges first in phase-space distance, while yielding also a small force residual, less than \textcolor{black}{$10\%$}. 
We also remark that the distance function has to be adapted in each case for optimized performance: the softer the material, the lower the $C$. 

%% file: subfiles/NeuralNetworks.tex
We also envision PSI solvers contributing to the adoption of neural-network based constitutive laws. 
In recent years, there has been a continuous stream \cite{Bessa:2019,Masi:2021,Mohr:2021,Masi:2022,Chen:2022,Rimoli:2022,review:2024} of ground-breaking work that has proved the capacity of deep learning to represent complex, history-dependent, non-linear material behavior. 
However, the computation of gradients with respect to the value of the entries (strain components) via ``automatic differentiation'' \cite{huot:2022} can be onerous and poorly-defined at times \cite{BECK:1994,Lee:2020,Lee:2023}, 
e.g., since neural networks tend to feature activation functions with limited continuity \cite{Apicella:2021} (for instance, ReLu, see \Cref{fig:approximation_NN}(g)). 
This logically puts in question the teaming up of neural-network constitutive laws with Newton-Raphson solvers. 

The following tests were carried out in a Python implementation of the same code used in the rest of the paper, but in a different laptop machine with the following specifications: 

\begin{itemize}[leftmargin=*]
    \item CPU: 12th Gen Intel(R) Core(TM) i5-1235U @ 4.40GHz, 10 cores.
    \item RAM: 16 Go, DDR4, 3200 MHz
    \item OS: Ubuntu 24.04.1 LTS.
\end{itemize}

The truss structure remains the same but the loading changes: the intensity of the imposed forces at the upper level, recall \Cref{fig:truss_bars}, is increased (see Supplementary material for details). 

\subsection{Usage and speed limitation of the neural-network}

\label{Sec:NN_limitations}

The neural-network used in this paper was trained on a pre-generated dataset. This dataset consisted of 1000 points distributed along a curve generated using \cref{eq:material_law}. In order to get a little closer to a real case based on experimental measurements, noise was added to all the points in this dataset, following the same strategy used by Kirchdoerfer and Ortiz \cite{Trent_1} (reproduced in \Cref{app:NeuralNetwork}).

\begin{figure}[h]
\centering
\captionsetup[subfigure]{justification=centering}
    \begin{subfigure}[t]{.49\linewidth}
    \includegraphics[width=\linewidth]{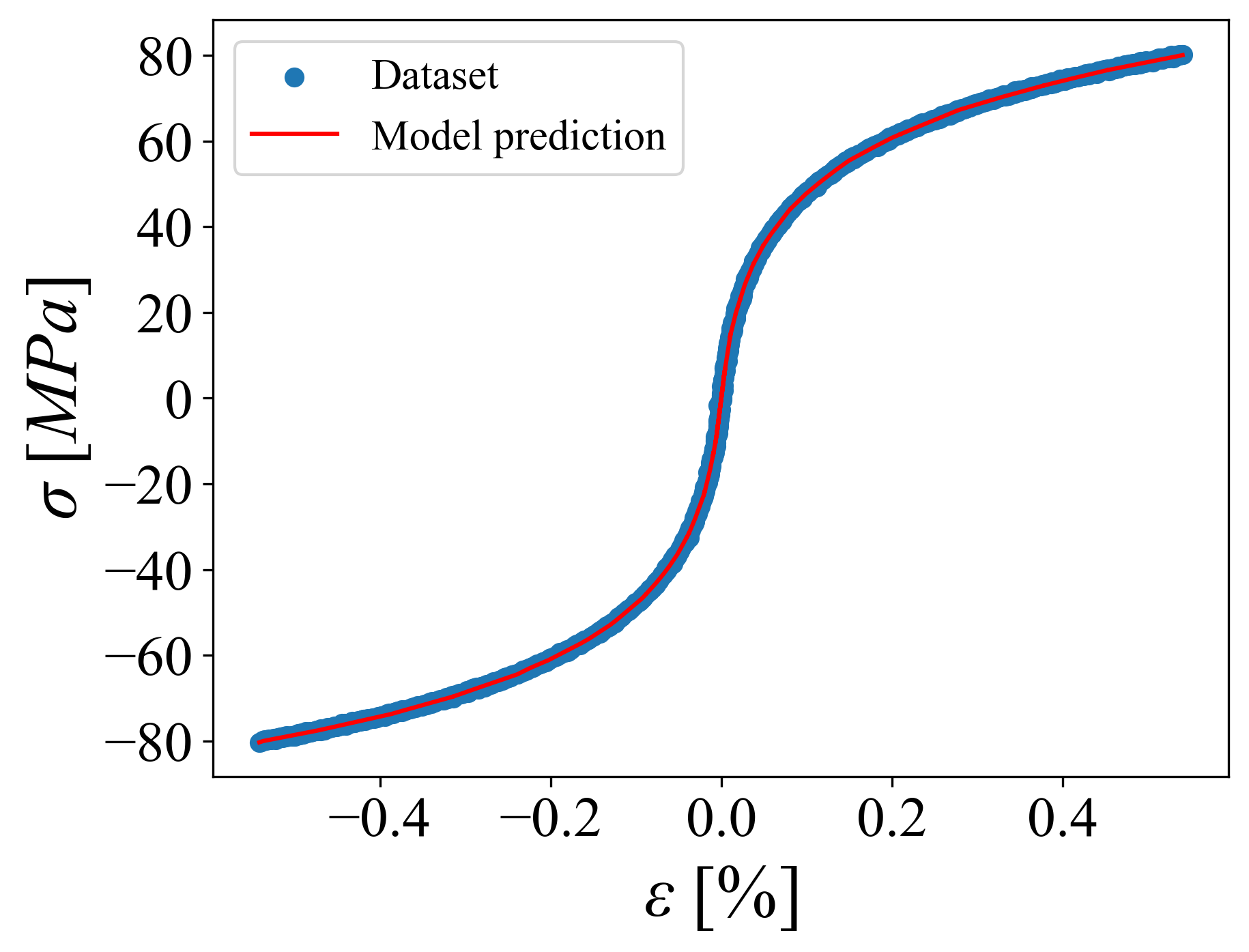}
        \caption{\Large (a)}
    \end{subfigure}
    \begin{subfigure}[t]{.49\linewidth}
    \includegraphics[width=\linewidth]{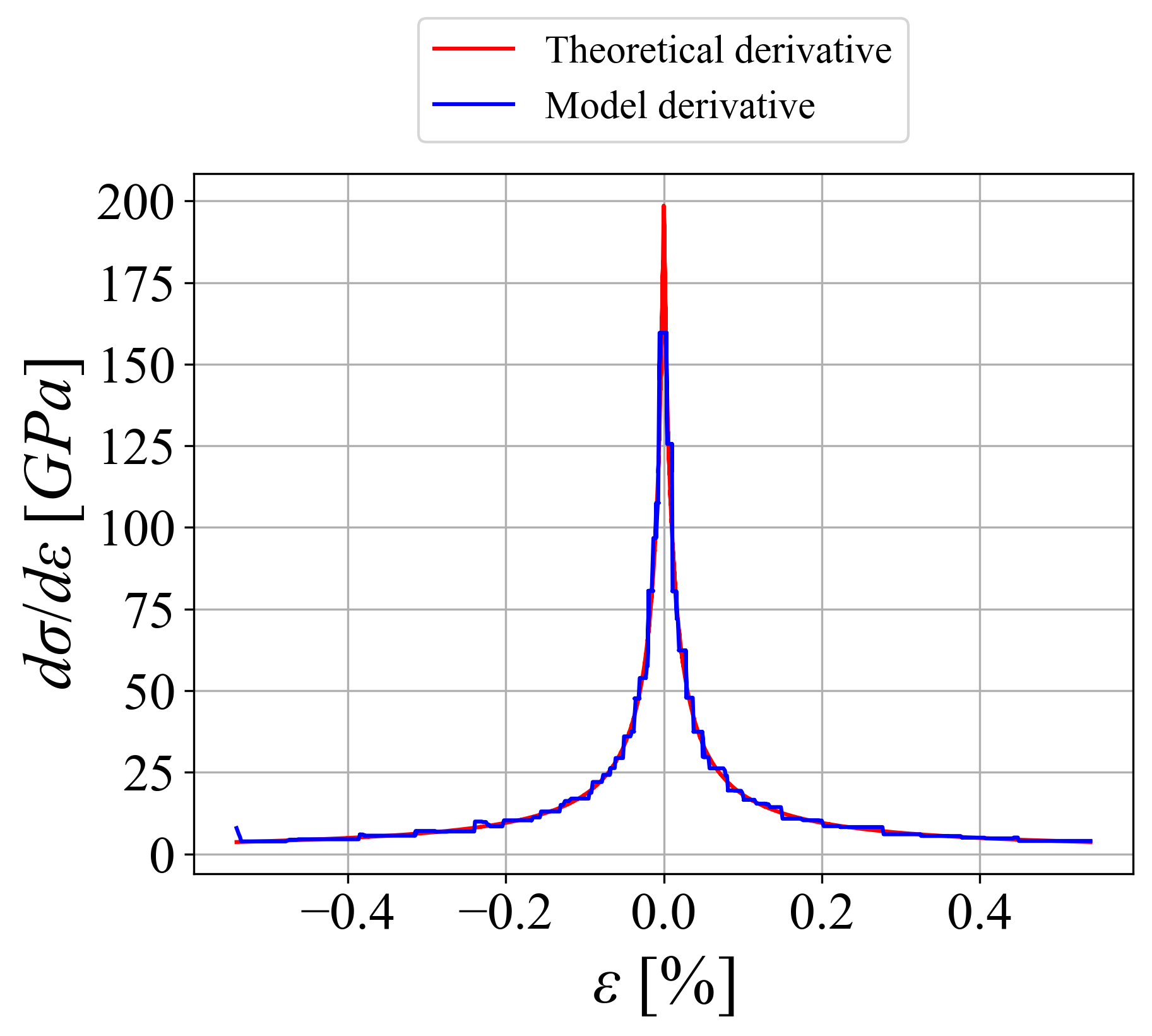}
        \caption{\Large (b)}
    \end{subfigure}
\caption{Estimation of the constitutive law using a neural network based on a noisy dataset, with (a) the representation of the model compared with the dataset (b) the derivative of the model compared with the derivative of the constitutive law.}
\label{fig:approximation_NN}
\end{figure}

Hyperparameter optimization using Optuna \cite{optuna_2019} and analysis using MLFlow \cite{mlflow} resulted in a neural network that best approximated the dataset as a continuous function. The neural network architecture consists of neurons activated with a ReLU (rectified linear unit) function. The strain is given as input, then distributed in three hidden layers, each comprising 112 fully connected neurons. Finally, the stress is recovered as output by passing through a neuron with a linear activation function. This gives a total of 25649 trainable parameters. The neural-network and its training is more detailed in \Cref{app:NeuralNetwork}. 
The final network result is depicted in \Cref{fig:approximation_NN}.

Because of the large number of elementary operations to be performed at every forward pass through the network, evaluation of the neural-network will be much slower and more costly than evaluating a customary material law (i.e. \cref{eq:material_law}). 
The evaluation time for the model derived from the neural-network is around $30\times$ longer than the time required to evaluate it using the material law. 
This results in a significant slowdown in both NR and PSI. NR has to calculate the model gradient for each element at each iteration. 
For PSI, this slowdown is even more noticeable: 
the model is recalled many times to solve the minimization problem used to project onto the materially-admissible set ($P_D$).

These intrinsic limitations of neural networks lead to a noticeable slowdown in the FEM resolution of the proposed system, but do not hinder its success. 
All the model evaluations are carried out directly using the Pytorch library, as are all the automatic differentiations. The projection onto the $P_D$ space is carried out with the help of the minimization function of Scipy library (\texttt{scipy.optimize.minimize}), using the conjugate gradient minimization method. 
%



\subsection{PSI-NR comparison}

As we have seen, the PSI method allows the FEM to be solved without ever having to calculate the gradient of the constitutive law, thus avoiding this problem. 
The downside, however, is that the model has to be called much more often to satisfy the algorithm. 
\Cref{tab:NN_PSI_vs_NR} compares the performance in terms of calculation time and number of iterations for the Newton-Raphson method and the PSI solver. 
Each simulation was run five times, and the result is the average of these five simulations.

\bgroup
\def\arraystretch{1.5}
\begin{table}[H]
\centering
\caption{Wallclock time comparisons and iterations using NR and PSI (1CPU) associated with a neural-network model.}
\begin{tabularx}{0.75\textwidth}{ c| *{3}{Y|} }
\cline{2-3}
 &  \multicolumn{1}{c|}{Newton-Raphson} & \multicolumn{1}{c|}{Phase-Space Iterative solver} \\ \hline
 \multicolumn{1}{|c|}{Running time [s]}& 32  & 119 \\ \hline
 \multicolumn{1}{|c|}{Iterations}& 84 & 57 \\ \hline
\end{tabularx}
\label{tab:NN_PSI_vs_NR}
\end{table}
\egroup

\begin{figure}[h]
\centering
\captionsetup[subfigure]{justification=centering}
    \begin{subfigure}[t]{.49\linewidth}
    \includegraphics[width=\linewidth]{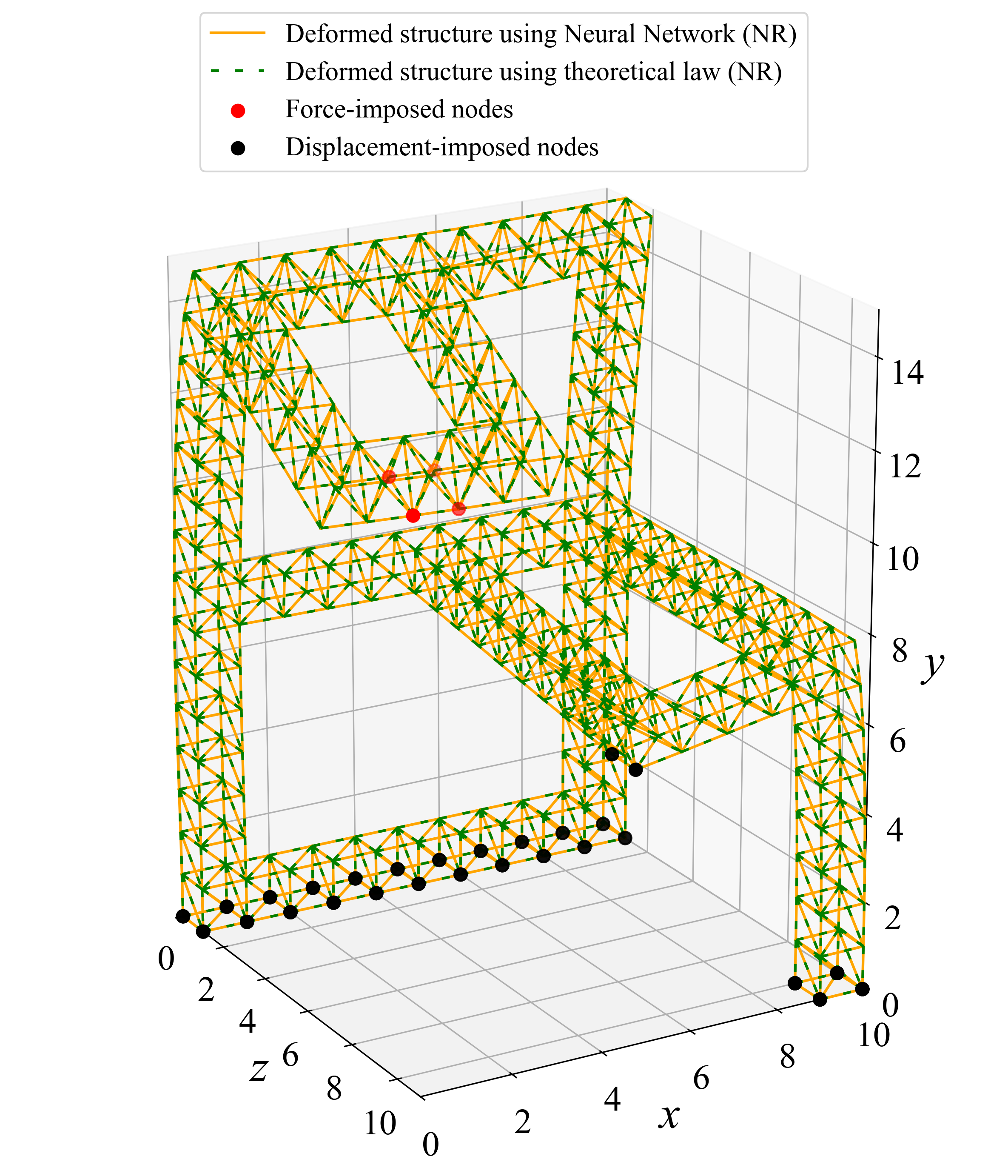}
        \caption{\Large (a)}
        \label{fig:Deformation_NN_NR}
    \end{subfigure}
    \begin{subfigure}[t]{.49\linewidth}
    \includegraphics[width=\linewidth]{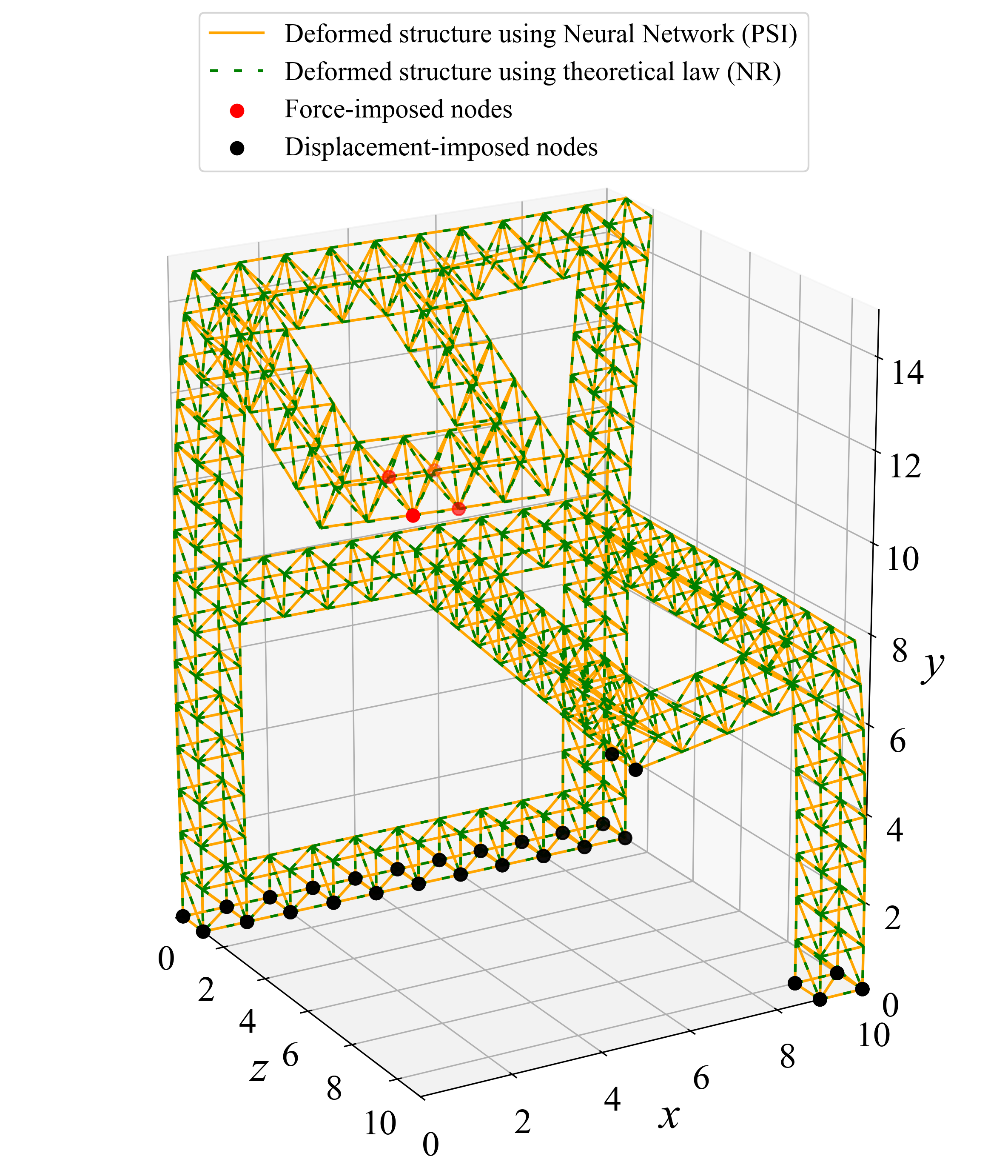}
        \caption{\Large (b)}
        \label{fig:Deformation_NN_PSI}
    \end{subfigure}
\caption{Units in meters. Deformation with a visual multiplier of x17, comparing the deformation obtained with Newton-Raphson applied from the constitutive law (with p = 0.0001) and (a) Newton-Raphson applied from the neural-network (b) a Phase-Space Iterative Solver applied from the neural-network}
\label{fig:Deformation_NN}
\end{figure}

The deformed structures using the different laws and methods are shown in Figure \ref{fig:Deformation_NN}. The forces applied to the red nodes are 1000 N, and the imposed displacement of the top black nodes is 0.05 m. The results, even when visually amplified, are virtually indistinguishable. 
%
%
%
The small differences
can be explained by the fact that the PSI stops once tol$_2$ is achieved (the change between each iteration is too small, see \Cref{sec:consecutive_iterations}). 
By decreasing the tolerance, better accuracy should be obtained, but at the cost of a longer calculation time. 
In the present case, tol$_1$ = 0.01 and tol$_2$ = 0.001. 
On the other hand, the error between the neural network and the constitutive law can be explained by two reasons: firstly, the fact that the neural network can never, by definition, approximate a law perfectly. 
The second reason is the introduction of noise into our dataset to simulate real measurements, generating an average standard deviation between the constitutive law and the dataset itself of $7.4\cdot 10^{-6}$ [\%] along $\varepsilon$ and $1.5 \cdot 10^{-4}$ [Pa] along $\sigma$.

As explained, the calculation time with PSI is longer, but can also be easily parallelized, due to the fact that the method is performed individually element by element at each iteration. We can therefore take advantage of using several cores of the CPU in parallel to reduce the calculation time. The results are shown in \Cref{tab:NN_PSI_CPU_paral}. All test were done 5 times, and averaged.


\bgroup
\def\arraystretch{1.5}
\begin{table}[H]
\centering
\caption{Wallclock time comparisons using PSI and multiple cores of the CPU associated with a neural-network model.}
\begin{tabularx}{0.75\textwidth}{ c| *{5}{Y|} }
\cline{2-5}
&  \multicolumn{4}{c|}{Phase-space iterative solver with NN}\\
\cline{2-5}
&  \multicolumn{1}{c|}{1 core} & \multicolumn{1}{c|}{2 cores} & \multicolumn{1}{c|}{4 cores} & \multicolumn{1}{c|}{8 cores}\\
\hline
 \multicolumn{1}{|c|}{Running time [s]}& 119 & 87 & 62 & 48\\ 
 \hline
\end{tabularx}
\label{tab:NN_PSI_CPU_paral}
\end{table}
\egroup

By using all 8 CPU cores, computing time is reduced by around 60\%. 
%
%



%



%% file: subfiles/Discussion.tex
Based on the examples analyzed above, there are three salient advantages in using PSI, all derived from the different treatment of the constitutive law: numerical handling of non-linearity, scalability, versatility. 
\begin{itemize}[leftmargin=*]
    \item Handling of non-linearity: since PSI does not update the stiffness matrix (nor its inverse), it avoids the poor-conditioning problem associated to local loss of stiffness in some elements \cite{Cook:2001,Oden:2013,Crisfield:2012}. 
    \item Scalability: whenever the $P_D$ projection is the performance bottleneck, we expect almost-linear time reduction when parallelizing this step and increasing the number of processors. 
    Moreover, since this step is carried out on a per-element basis, the gain should scale linearly with mesh size (number of elements). 
    \item Versatility: the two projections are but two minimizations, there is a bevy of numerical methods available to solve them. 
    For the projection onto the physically-admissible set, enough continuity can be always assumed so that the corresponding Euler-Lagrange equations can be found, what yield a simple linear system in all cases. 
    On the contrary, for non-linear constitutive laws, the per-element projection onto the materially-admissible set is intrinsically non-linear, so the user can opt among different ways of proceeding. 
    Depending, on the dimensionality of the local phase space or the available libraries' performance \cite{petsc-user-ref}, one could choose among multiple ways to carry out the projection onto the constitutive law. 
    When it comes to neural-network constitutive laws of material response, PSI scores yet another point given that it can work directly on the strain-stress pairs, without ever requiring derivatives of the constitutive law, which could be poorly computed even when using automatic differentiation (\Cref{fig:approximation_NN}).  
    %
    %
\end{itemize}

On top of everything said above, the mathematical underpinnings of the method can still be better understood if we extend the insights discussed in \Cref{sec:mathematical_props}. 
We also foresee borrowing ideas from POCS that are known to work in improving convergence. 
For instance: the theory of alternating projections goes beyond the simple two-projection setting we have been using \cite{Kopecka:2014}, 
some preliminary tests point to adding an intermediate third projection, consisting simply of evaluating the constitutive law, as a direct means to speed convergence up.    



We also point out the possibility of splitting the system in two regions: one that remains linear-elastic and another one where inelastic behavior manifests. 
A classic linear-elastic solver for the former can be combined with phase-space iterations in the latter \cite{Method_1}; this could lead to even faster simulations of that kind of systems where non-linear material behavior is localized in space. 
Going one more step forward to an incremental-load scenario, we can envision localizing non-linearity both in space and time, where the parts of the domain that are solved in the traditional way or with phase-space iterations evolve dynamically with the load \cite{Wattel:2023}.    


%% file: subfiles/Remarks_Acknow.tex

We have introduced the concept of phase-space solvers, a new kind of iterative algorithm for solving problems in elasticity, and have shown its potential to outperform traditional solvers when it comes to solve large systems with severe non-linearity. 
The method relies on projecting onto convex sets whose frontier is defined either by the material constitutive law or Newton's second law. 
While the latter is achieved through solving a linear system of equations, the projection onto the constitutive law demands the element-wise resolution of either a non-linear optimization problem or a non-linear algebraic equation, but this step can be trivially parallelized, 
what renders the algorithm fast enough to appear competitive against other options. 
See that we have not even started to exploit the advantages of GPUs for this purpose \cite{GPU}, but this certainly is in the horizon, among other improvements. 

Further substantiating the PSI strengths that have been discussed already
is our next goal, hence we are to implement PSI in a dedicated, well-maintained, thoroughly-tested solid mechanics package \cite{Akantu,Deal:II,FENICS:book} in order to carry out comprehensive performance benchmarking. 
It remains to formally characterize errors \cite{Kayalar:1988} and convergence rates \cite{bauschke:2003,badea:2012,BADEA:2016} in more general settings (e.g., large meshes, using different distances, element types, adding extra projections), but this seems a feasible task since much inspiration can be drawn from previous work on the method of alternating projections or the method of projections onto convex sets, and preliminary results have already been laid out in this paper. 

Extra physics to be considered include dynamics \cite{Trent_2} and finite kinematics \cite{Oden:2013,Zienkiewicz:1991}, apart from inelastic (irreversible) problems with complex non-monotonous loading paths or material instabilities. 
For the latter, we foresee synergies with novel non-convex optimization techniques \cite{Goodfellow:2016}. 
Finally, we highlight that this approach can also be used to solve non-linear problems in other branches of physics beyond mechanics, as long as they feature a constitutive law and some physical balance conditions (Newton's second law, conservation of electric charge...)  
As an example, let us mention the analysis of electric circuits \cite{Galetzka:2024}. 

\section*{Supplementary material}

The Mathematica notebooks used to carry out simulations in \Cref{sec:truss} can be retrieved from the repository named \texttt{PSI} in the last author's GitHub page \texttt{github.com/jgarciasuarez}.

The Jupyter notebooks and Python script used for \cref{sec:NN_gaetan} can be retrieved from \texttt{gitlab.epfl.ch/gcortes} \texttt{/nn\_iterative\_solver}.

\section*{Acknowledgments}

The last author is thankful to Dr. Trenton Kirchdoerfer for sharing information about the homonymous truss. 
We would like to thank also Prof. Laurent Stainier, who brought the LATIN method to our attention. 
The comments of Prof. Torres de Lizaur, regarding the scope of Section 3, have been instrumental and are gratefully acknowledged. 
We also note useful discussions with Prof. J.-F. Molinari, Dr. Lucas Fourel, Dr. Guillaume Anciaux and Dr. Nicholas Richart. 
The work of the last and first authors is sponsored by the Swiss National Science Foundation via the Ambizione Grant 216341 ``Data-Driven Computational Friction'', awarded to the last author, for which we are exceedingly grateful. 

%% file: subfiles/appendix/Maths.tex

\begin{definition}[Scalar Product]
An inner or scalar product, in a vector space $V$ over the field $F$, is, a map $\langle \cdot,\cdot \rangle : V \times V \rightarrow F $ satisfying the following properties $\forall x,y,z \in V$ and $\forall$ scalars $a,b \in F$:
\begin{itemize}
    \item symmetry: $\langle x,y \rangle$=$\overline{\langle y,x \rangle}$. If $F$ is $\mathbb{R}$, conjugate symmetry is just symmetry.
    \item linearity: $\langle ax+by,z\rangle$=$a\langle x,z\rangle$+$b\langle y,z\rangle$
    \item positive definiteness: if $x$ is not zero, then $\langle x,x \rangle >0$
\end{itemize}
    
\end{definition}


\begin{definition}[Norm]

Given a vector space $\mathcal{X}$ over a subfield $F$ of the complex numbers, a norm of $\mathcal{X}$ is a rel function p: $X \rightarrow \mathbb{R}$ with the following properties:
\begin{itemize}
    \item Triangle inequality: $p(x+y) \leq p(x)+p(y)$, 
       $\forall x,y$ $in$ $\mathcal{X}$
     \item Homogeneity: $p(sx)=|s|p(x), \forall$ $x$ $in$ $\mathcal{X}$ and scalar $p$
     \item Positive definiteness: $\forall$ $x$ $in$ $\mathcal{X}$, $x=0$ if and only if $p(x)=0$
\end{itemize}

\label{def:norm}
    
\end{definition}

\begin{definition}[Hilbert space]
A Hilbert space is a real or complex inner product space that is also a complete metric space with respect to the distance function induced by the inner product.
\label{def:hilbert}
\end{definition}


%% file: subfiles/appendix/NR.tex
As most iterative solvers, this procedure aims at minimizing the residual $r = \norm{\Delta F(\boldsymbol{u})} = \norm{\boldsymbol{F}_{int}(\boldsymbol{u}) - \boldsymbol{F}_{ext}}$. 
This is done in incremental fashion, starting from some initial guess, e.g., $\boldsymbol{u}_0 = \boldsymbol{0} $: assuming the solution at the $n$-th loading step is known ($\boldsymbol{u}_n $), then 
\begin{align}
    \Delta \boldsymbol{u}_{n} 
    = 
    - \boldsymbol{T}(\boldsymbol{u}_n)^{-1} \Delta \boldsymbol{F}_{n} 
    = 
    - \boldsymbol{T}(\boldsymbol{u}_n)^{-1} 
    \left[ \boldsymbol{F}_{int}(\boldsymbol{u}_n)  - \boldsymbol{F}_{ext} \right] \, ,
    \label{eq:NR_update}
\end{align}
where $\boldsymbol{T}(\boldsymbol{u}_n)$ is the tangent stiffness computed from the prior displacement field, and this $\boldsymbol{u}_{n+1} = \boldsymbol{u}_{n} + \Delta \boldsymbol{u}_{n}$. The tangent stiffness matrix is, in analytical form, 
\begin{align}
    \boldsymbol{T}(\boldsymbol{u}_n)
    =
    {\partial \boldsymbol{F} 
    \over 
    \partial \boldsymbol{u}}
    \Big|_{\boldsymbol{u} = \boldsymbol{u}_n} \, ,
\end{align}
which, in practice, must be re-assembled 
at every iteration and involves derivatives of the constitutive law, computed either analytically or numerically. 
 
Oftentimes, for efficiency purposes, the new tangent stiffness is weighted with the one from prior steps or with the ``undeformed'' stiffness, i.e., $\boldsymbol{T}(\boldsymbol{u}_n) \to \gamma \boldsymbol{T}(\boldsymbol{u}_n) + (1 - \gamma) \boldsymbol{T}(\boldsymbol{0})$.

For the sake of completeness, let us also mention the quasi-NR method, which computes the inverse tangent matrix only for the first iteration and then uses Sherman-Morrison formula to update this inverse, assuming rank-one matrices. 
%

A word on the enforcement of Dirichlet BCs in the implementation of NR used in this work. 
We have opted for an approach that sets the known values directly on the initial guess, while subsequent iterations solved the incremental system \cref{eq:NR_update} only for the free dofs. 
This is a popular procedure, as its implementation is straighforward, but it can challenge the convergence of the method: 
having a strongly localized displacement field in the first guess around the nodes with essential BCs may lead to a deficient tangent stiffness matrix. 
This, in turn, may lead to nonsensical solutions or to difficulties when it comes to convergence. 
An antidote to this issue is precisely weighting tangent stiffness at every step with the zero-deformation stiffness matrix $\boldsymbol{T}(\boldsymbol{0})$. 

%% file: subfiles/appendix/NeuralNetwork.tex
\setcounter{figure}{0}
\renewcommand{\thefigure}{A.\arabic{figure}}

The aim of this chapter is to describe the process of designing and training the neural network used in this paper. 

\subsection{Generating the dataset}

Firstly, a training dataset was generated from sampling a given function. It is defined using the equation below:

\begin{equation}
    \sigma = m(\epsilon) = E_0((|\epsilon|+c))^p-c^p)*\mathrm{sign}(\epsilon)
    \label{eq:material_law_NN}
\end{equation}

\begin{figure}[h]
    \centering
    \includegraphics[width=0.5\linewidth]{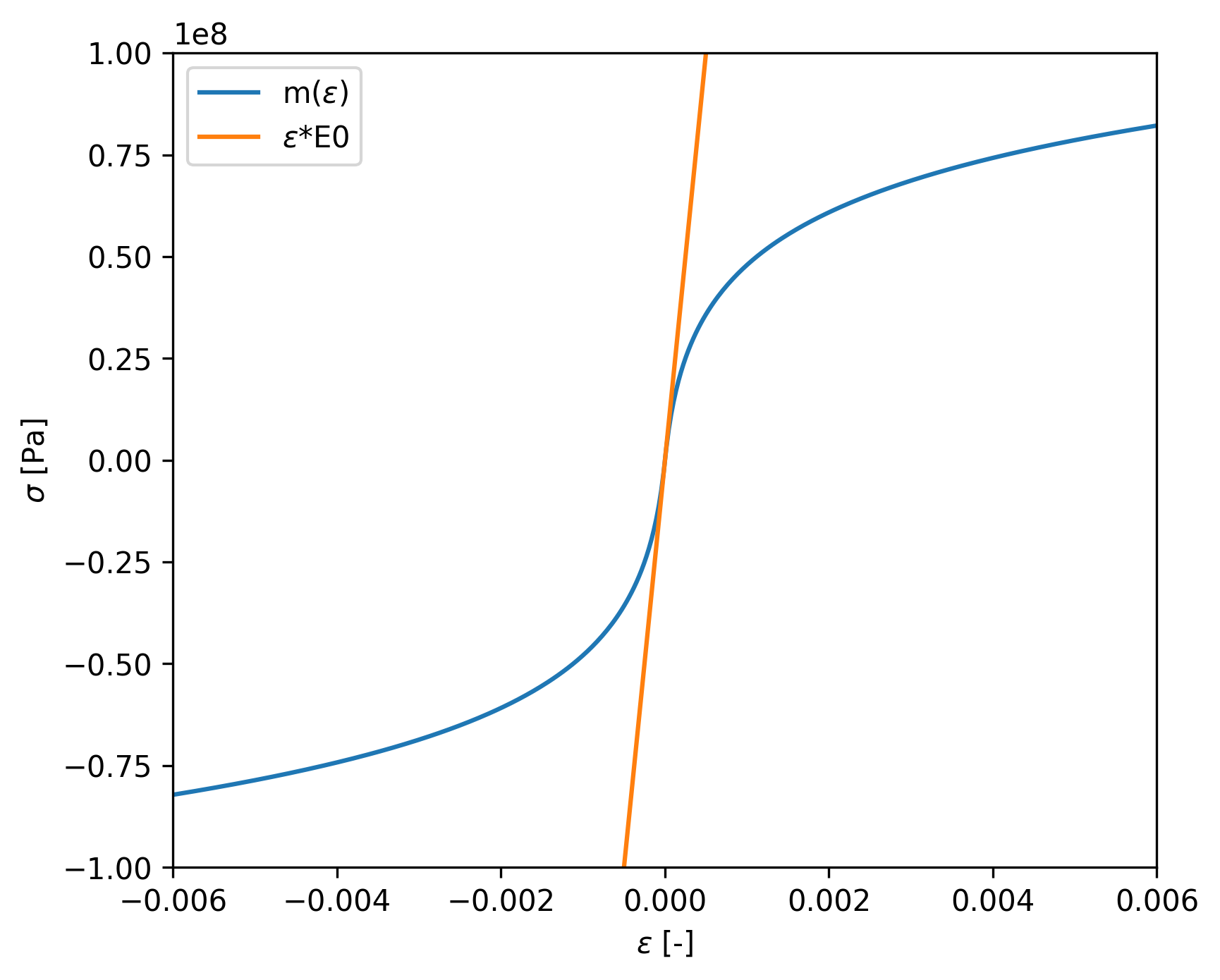}
    \caption{Material law defined by equation \cref{eq:material_law_NN}. In blue, the constitutive law. In orange, the slope of the constitutive law at $\epsilon = 0$. }
    \label{fig:material_law_NN}
\end{figure}

With $p = 0.0001$, $c = \frac{1}{p}^\frac{1}{p-1} = p^{(1-p)^{-1}} = 9.991$ and $E_0 = 2 \cdot 10^{11}$ Pa. The material law is shown in \Cref{fig:material_law_NN}. 
To create the dataset, the generation interval was first determined. To do this, the results of the NR algorithm based on the constitutive law were studied (\Cref{sec:truss}): 
the maximum stress experienced throughout the structure defined the upper and lower limit required for our dataset. By multiplying it by 1.2 to allow a little margin (particularly during the iterative phase), we arrive at a dataset generated between $\epsilon =$ -0.0054 and $\epsilon =$ 0.0054.

The dataset should consist of enough points to train the neural network. 
To choose where to place, say, 1000 points on the curve, the simplest solution would have been to generate 1000 points at regular intervals along $\epsilon$. 
However, this solution can be problematic due to the large slope close to 0, which makes it difficult for our neural network to approximate the function around the origin. 
Instead, 1000 points should be created at regular intervals along the arc, so that the points in the dataset are always at equal distances from each other.

In order to make this dataset more realistic, noise has been added. 
There are several possible solutions in a case like this: the noise can be proportional, constant, or relative to the total number of points. 
The latter solution was chosen, as it is similar to the way in which Kirchdoerfer et al.\cite{Trent_1} initially generated their dataset. 
To do this, the value $C_e$ is defined to readjust the scale of the constraint $\sigma$ and have it in the same order of magnitude as $\epsilon$: $C_e = \frac{E_0}{10}$. This value is taken from the paper, and works well for having close values on the $\epsilon$ and $\frac{\sigma}{C_e}$ scale. 
This readjustment allows us to create noise proportional to the space between each point without scale problem, but also to readjust all the datapoints so that they are equidistant in the new space. 
Stress will therefore be just as important as strain in the distribution of points along the curve.

\begin{figure}[h]
\centering
\captionsetup[subfigure]{justification=centering}
    \begin{subfigure}[t]{.49\linewidth}
    \includegraphics[width=\linewidth]{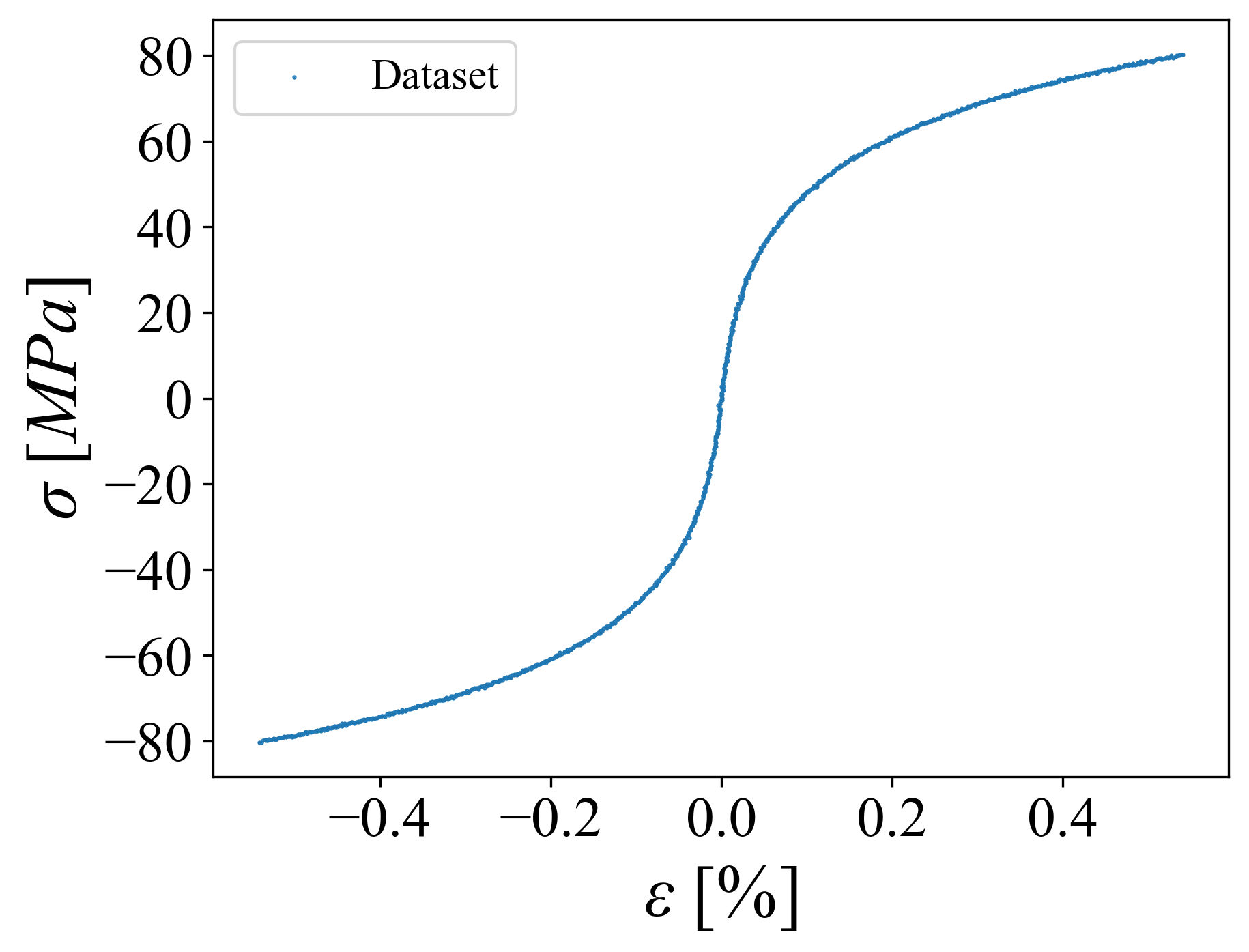}
        \caption{\Large (a)}
        \label{fig:dataset_1x}
    \end{subfigure}
    \begin{subfigure}[t]{.49\linewidth}
    \includegraphics[width=\linewidth]{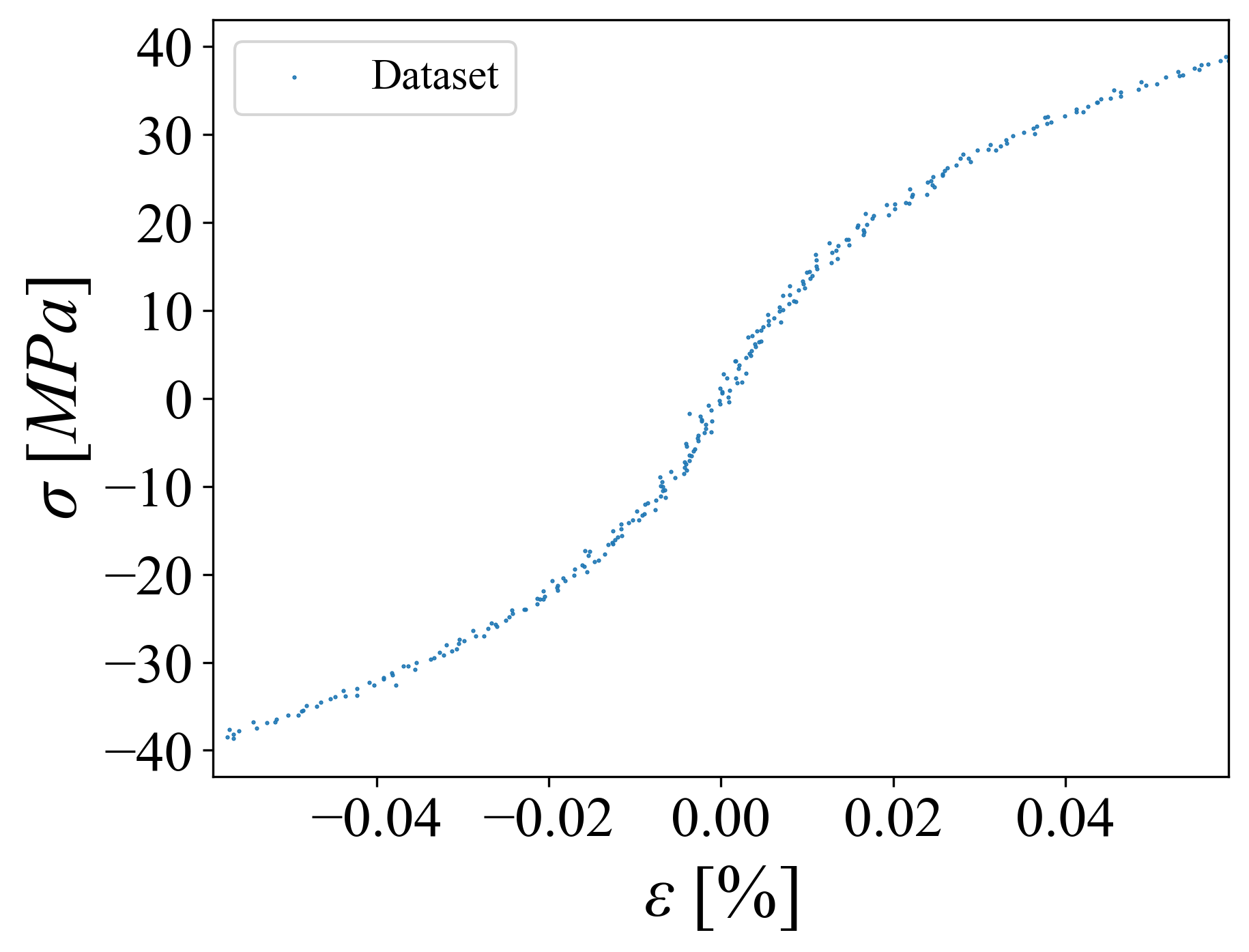}
        \caption{\Large (b)}
        \label{fig:dataset_10x}
    \end{subfigure}
\caption{Noisy dataset composed of 1000 points generated to train the neural network, showing (a) the whole dataset (b) a focus around the origin}
\label{fig:Dataset}
\end{figure}

The distance is calculated between each point on the curve in this new space, and noise is added to the set of points as a function of this distance. 
The calculation of the noise added to the points is based on a normal distribution with a standard deviation equal to the distance between each point on the curve, divided by two. 
Back to the original space, this gives an average standard deviation between the constitutive law and the dataset itself of $7.4\cdot 10^{-6}$ along $\epsilon$ and $1.5\cdot 10^{-4}$ Pa along $\sigma$. The dataset is shown in \Cref{fig:Dataset}.

\subsection{Creating and training the network}

Firstly, the hyperparameters and the architecture of the neural network had to be established. 
The primary aim of the neural network was to transform an input ($\epsilon$) into an output ($\sigma$), taking into account a strong non-linearity away from the origin. 
For this reason, a fully connected dense network was chosen, with a ReLU activation function. 
After preliminary tests, the maximum epoch limit of $10^4$ was deemed to be sufficient to achieve a good fit, with a patience parameter of $10^3$ epochs \cite{Goodfellow:2016}. 
Three more hyperparameters remained to be determined: the number of hidden layers, the number of neurons per layer, and the learning rate. 
These three were determined using Optuna \cite{optuna_2019}, by varying the number of hidden layers between 1 and 5, the number of neurons per layer between 30 and 120, and the learning rate between 10$^{-5}$ and 10$^{-3}$.

For the loss function, a simple MSE function was implemented, comparing the data output by the neural network with that of the dataset. 
Adam optimizer algorithm was used to train the neural network. 
The final value used to assess the neural network's effectiveness in processing the dataset is that of the loss function applied to a validation dataset. 
It is made up of 20\% of all the points in the dataset, randomly selected along the curve, leaving the remaining 80\% for training the neural network (i.e. 200 and 800 points). 
Finally, before being integrated into the training of the neural network, the entire dataset is normalised between 0 and 1. 
This normalization makes it easier and faster to train the neural network, avoiding too rapid a saturation of the neurons.

After analyzing the results in Optuna and MLFlow \cite{mlflow}, the combination of parameters that gave the lowest loss function was 3 layers, 112 neurons per layer and a learning rate of $6 \cdot 10^{-5}$. 

\begin{figure}[h]
    \centering
    \includegraphics[width=0.5\linewidth]{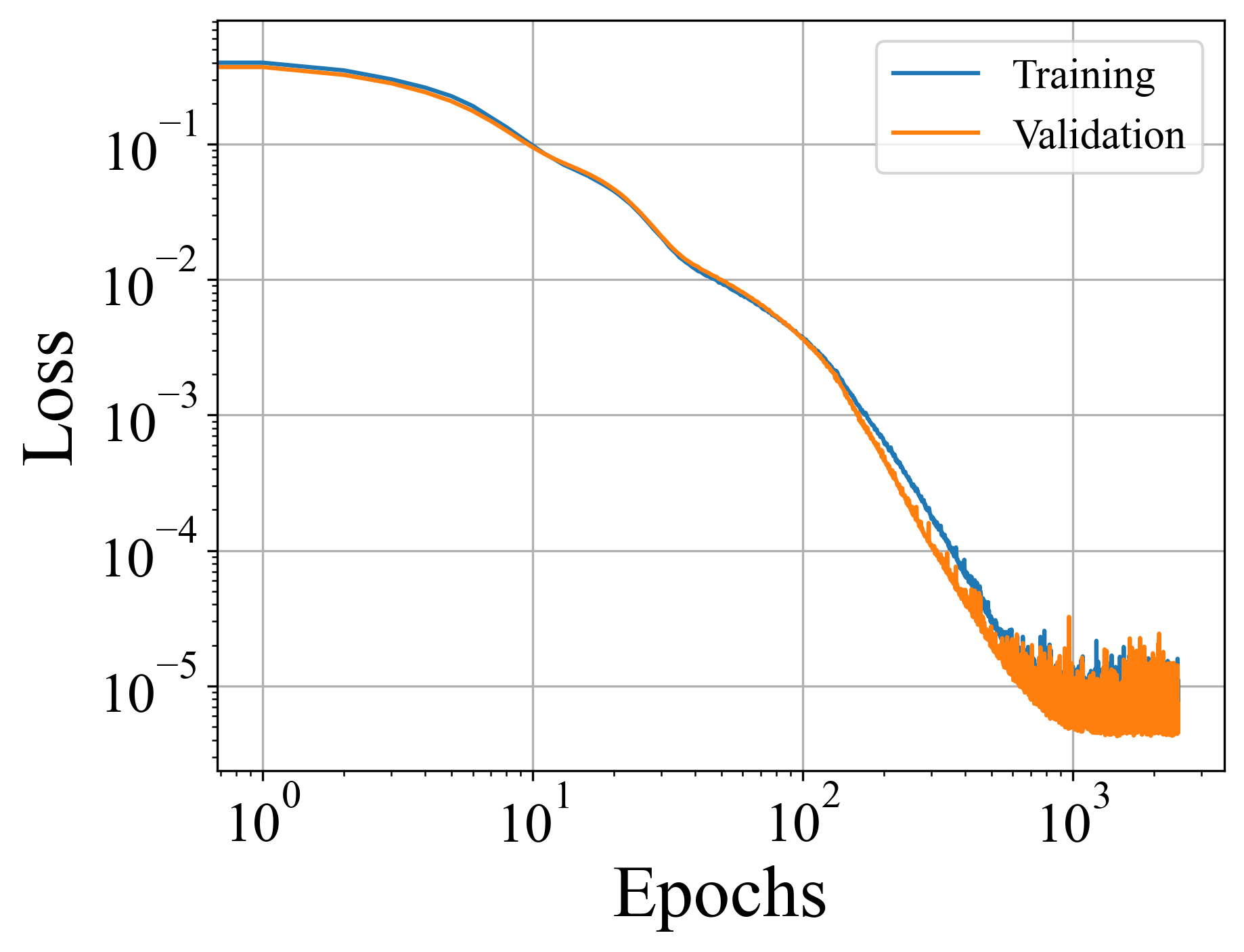}
    \caption{Evolution of the loss function as a function of the number of epochs. The loss function is calculated for the training set (in blue) and the validation set (in orange)}
    \label{fig:loss_function_NN}
\end{figure}

The evolution of the loss function as a function of time is shown in \Cref{fig:loss_function_NN}, and the representation of the model obtained in comparison with the dataset and the constitutive law is shown in \Cref{fig:approximation_NN}.

The loss function begins to plateau at around 1000 epochs, with the value for training and validation fairly equivalent, indicating no overfit \cite{Goodfellow:2016}. 
Due to the nature of the dataset, which is noisy, it is difficult to go any lower. 
These results are therefore satisfactory from a training point of view. 
As for the fit itself, the curve follows the dataset trend with good accuracy. 
However, if we analyse the derivative of the curve in more detail, we notice that it is very difficult to reproduce the peak around the origin, i.e., the zero-strain stiffness. 
This difficulty can be explained by several factors: firstly, the high non-linearity of the dataset leading to a very rapid increase in the slope close to 0. 
The ReLU functions will have difficulty adapting to this, but a more consistent dataset close to the highest slope could help to counter this problem. 
Another factor is noise: the steep slope may require fairly accur

Even if it does not reproduce the constitutive law perfectly, the result obtained is atisfactory in the context of a project where only data is available, with no constitutive law for comparison. 
The main drawback, however, is the high number of parameters (25649). This considerably increases the calculation time for each model evaluation. This problem is discussed in more detail in \Cref{Sec:NN_limitations}.

\subsection{Training without noise}

\begin{figure}[h]
\centering
\captionsetup[subfigure]{justification=centering}
    \begin{subfigure}[t]{.45\linewidth}
    \includegraphics[width=\linewidth]{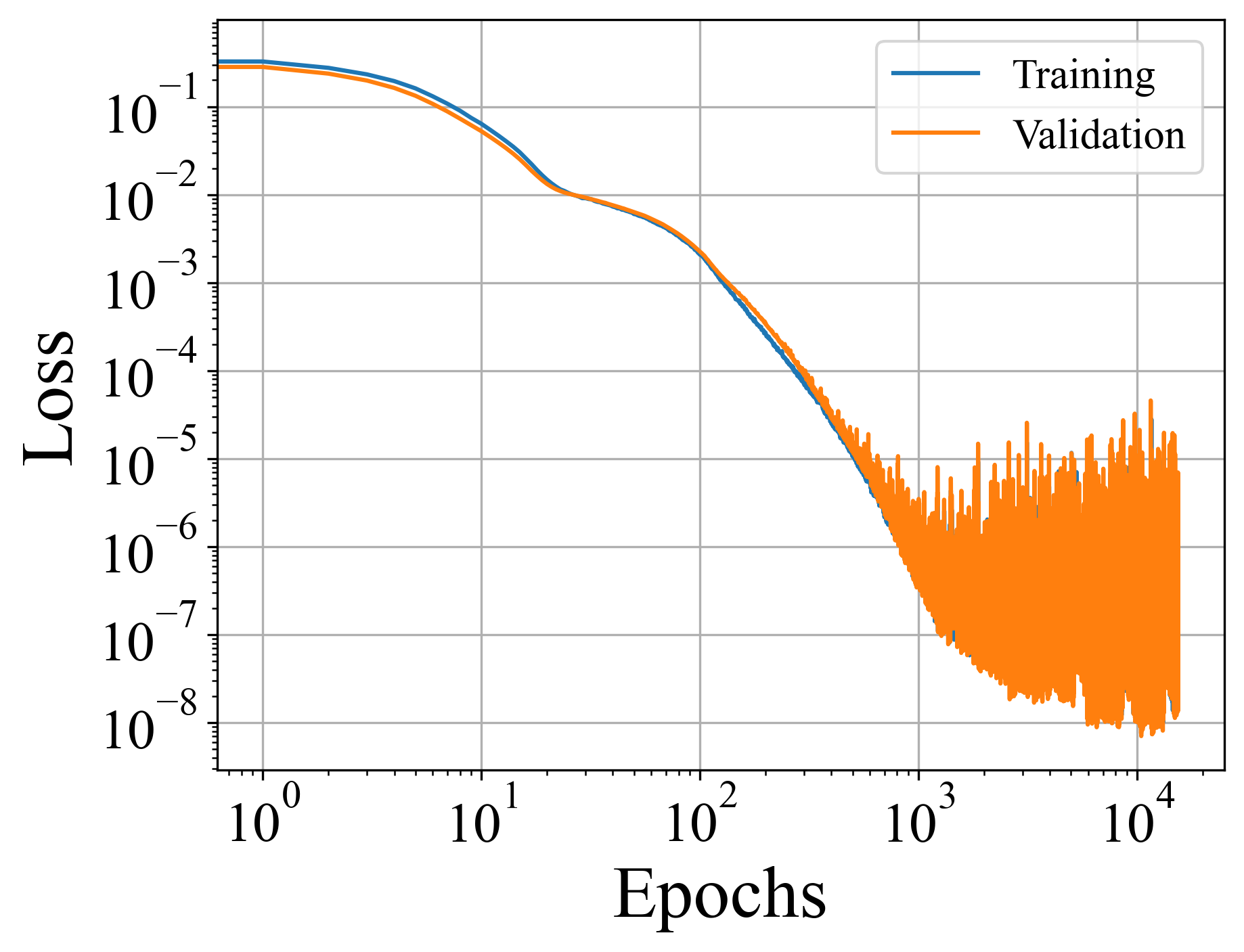}
        \caption{\Large (a)}
    \end{subfigure}
    \begin{subfigure}[t]{.45\linewidth}
    \includegraphics[width=\linewidth]{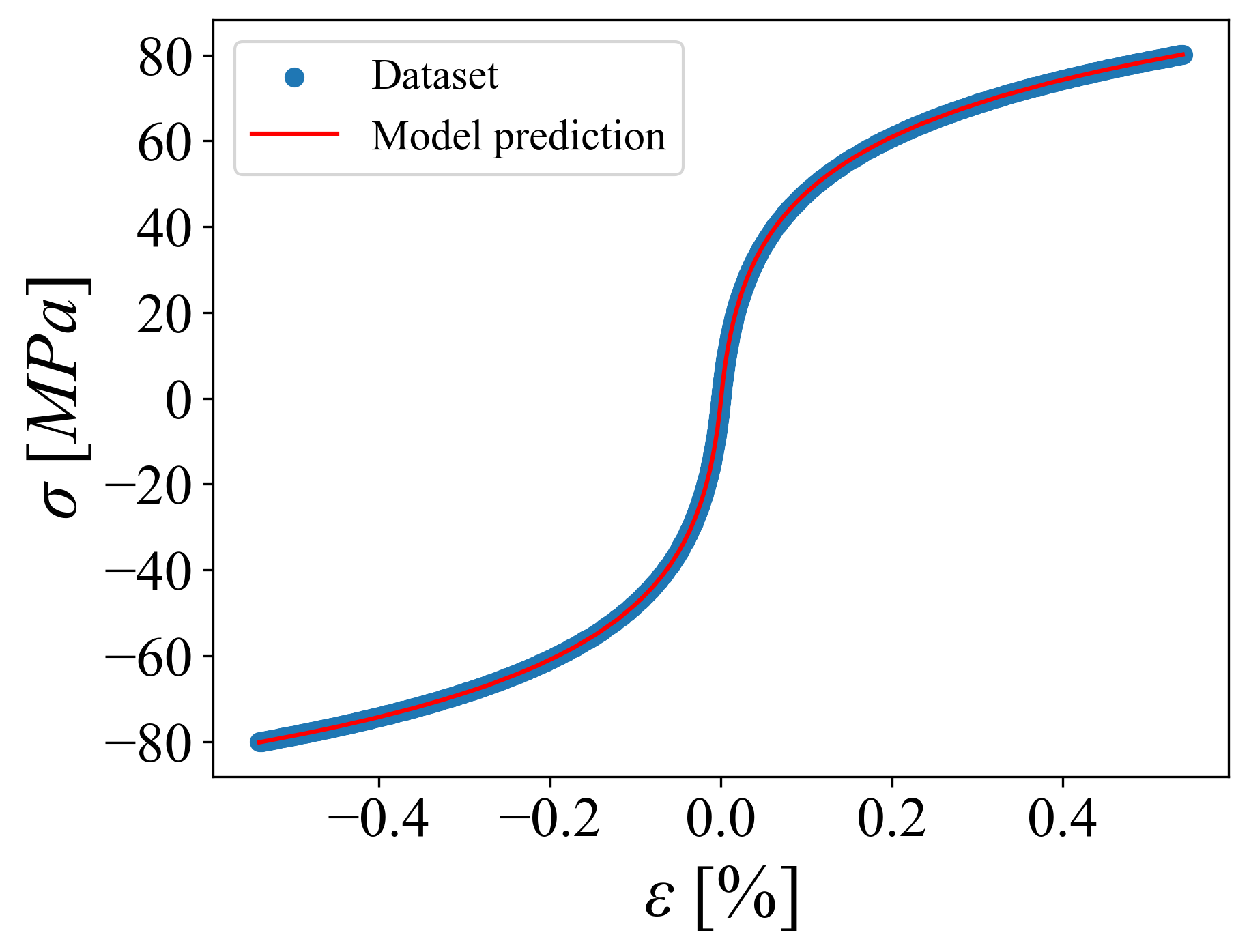}
        \caption{\Large (b)}
    \end{subfigure}
    \\
    \begin{subfigure}[t]{.45\linewidth}
    \includegraphics[width=\linewidth]{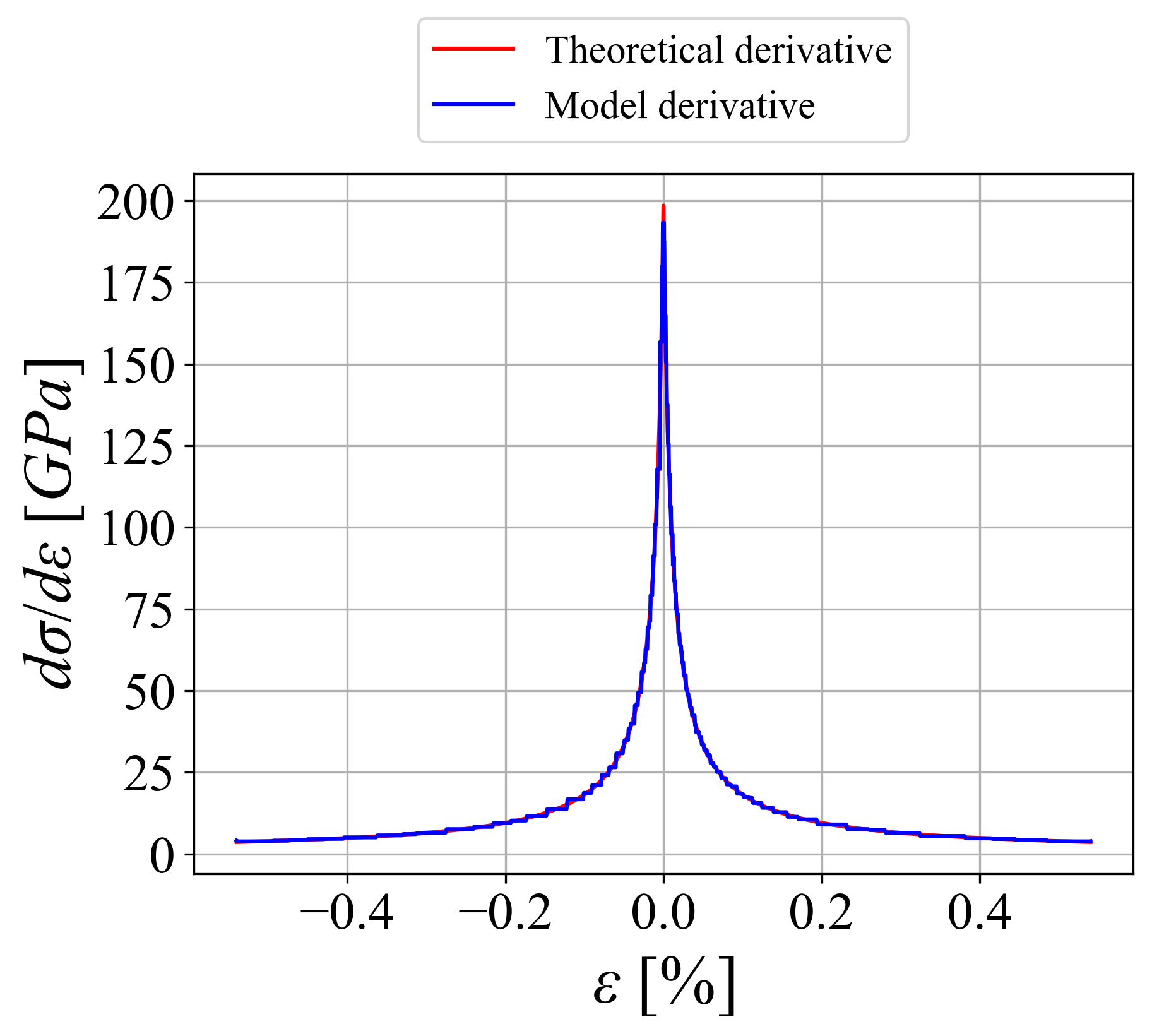}
        \caption{\Large (c)}
    \end{subfigure}
\caption{Results of training a neural network with a noise-free dataset, with (a) the evolution of the loss function (b) the approximation of the constitutive law (c) the derivative of the model.}
\label{fig:NN_noiseless}
\end{figure}

Finally, in order to gain a better appreciation of the neural network and better understand the difficulty of tracking the slope close to zero, a new dataset is generated.
It is composed of the same number of points (1000), this time no noise is generated, perfectly sampling the material's constitutive law.
For a better comparison, the same hyperparameters and architecture are used to train the neural network.
The only difference is the increase in patience parameter, from 1000 to 5000 epochs, because the loss function has a more noticeable continuous decrease.
The results obtained for the loss function, the curve approximation and the derivative are shown in \Cref{fig:NN_noiseless}.

We can see that the slope is almost perfectly tracked around zero, which tends to confirm the hypothesis that the difficulty in obtaining the right slope around zero with the original dataset is mainly due to noise.